\documentclass[11pt]{article}
\usepackage{amsmath,amssymb,amsbsy}
\newcommand{\beq}{\begin{eqnarray}}
\newcommand{\eeq}{\end{eqnarray}}
\newcommand{\bei}{\begin{itemize}}
\newcommand{\eei}{\end{itemize}}
\newcommand{\bee}{\begin{enumerate}}
\newcommand{\eee}{\end{enumerate}}
\newcommand{\beqe}{\begin{eqnarray*}}
\newcommand{\eeqe}{\end{eqnarray*}}

\newtheorem{prop}{Proposition}
\newtheorem{lemma}{Lemma}

\newtheorem{theo}{Theorem}
\newtheorem{cor}{Corollary}
\newcommand{\un} {\mbox{ \rm 1\hspace{-0.3em}I}}
\newcommand{\pa}[1]{\left({#1}\right)}
\newcommand{\cro}[1]{\left[{#1}\right]}

\newcommand{\ac}[1]{\left\{{#1}\right\}}

\newcommand{\EE}{\mathbb{E}}
\newcommand{\PP}{\mathbb{P}}
\newcommand{\F}{\mathcal{F}}
\newcommand{\hT}{\hat{T}}
\newcommand{\hm}{\hat{m}}
\newcommand{\M}{\mathcal{M}}
\newcommand{\N}{\mathbb{N}}
\newcommand{\cS}{\mathcal{S}}

\begin{document}
\renewcommand{\theequation}{\thesection.\arabic{equation}}
\title{ Non asymptotic minimax rates of testing in signal detection with heterogeneous variances}
\author{ B. Laurent, J.M. Loubes, C. Marteau \thanks{Institut de Math\'ematiques de Toulouse, INSA de Toulouse, Universit\'e de Toulouse.} }
\date{}
\maketitle
\begin{abstract}
The aim of this paper is to establish non-asymptotic minimax rates   for goodness-of-fit hypotheses testing in an heteroscedastic setting. More precisely, we deal with sequences $(Y_j)_{j\in J}$ of independent Gaussian random variables, having mean $(\theta_j)_{j\in J}$ and variance $(\sigma_j)_{j\in J}$. The set $J$ will be either finite or countable. In particular, such a model covers the inverse problem setting where few results in test theory have been obtained. The rates of testing are obtained with respect to $l_2$ and $l_{\infty}$ norms, without assumption on $(\sigma_j)_{j\in J}$ and on several functions spaces. Our point of view is entirely non-asymptotic.
\end{abstract}


\section{Introduction}
We consider the following heteroscedastic statistical model :
\begin{equation}
\label{mod1}
Y_j=\theta_j+\sigma_j \epsilon_j,\quad j \in J,
\end{equation}
where $\theta =( \theta_j)_{j \in J} $ is unknown, $\sigma =( \sigma_j)_{j \in J} $ is assumed to be known, and the variables $(\epsilon_j)_{j \in J} $ are i.i.d. standard normal variables. The set $J$ is either $\{1, \ldots,N\} $ for some $N \in \mathbb{N}^*$ (which corresponds to a Gaussian regression model) or $\mathbb{N}^*$ (which corresponds to the Gaussian sequence model). The sequence $\theta$ has to be tested from the observations $(Y_j)_{j\in J}$ in order to decide whether $"\theta=0"$ or not. The particular case $\sigma_j=\sigma$ for all $j\in J$ corresponds to the classical statistical model where the variance of the observations is always the same. It has been widely considered in the literature, both for test and estimation approaches. In this paper, we consider a slightly different setting in the sense that the variance of the sequence is allowed to depend on $j$.\\ \indent We point out that  the model (\ref{mod1}) can describe inverse problems. Indeed, for a linear operator $T$ on an Hilbert space $H$ with inner product $(.,.)$, consider an unknown function $f$ indirectly observed in a Gaussian white noise model 
\begin{equation}
\label{eq:invfonc} 
Y(g)=(Tf,g)+ \sigma \epsilon(g), \: g \in H 
\end{equation} 
where $\epsilon(g)$ is a  centered Gaussian variable with variance $\|g\|^2:=(g,g)$. If $T$ is assumed to be compact, it admits a singular value decomposition (SVD) $(b_j,\psi_j,\phi_j)_{j \geq 1}$ in the sense that
$$T\phi_j=b_j \psi_j,\quad T^\ast \psi_j=b_j \phi_j,$$ with $T^\ast$ the adjoint operator of $T$. Hence considering the observations $Y(\psi_j)$, model (\ref{eq:invfonc}) becomes 
\begin{equation}
Z_j=b_j\theta_j +\sigma  \epsilon_j, \ j \in \N^{\star},
\label{mod2}
\end{equation} 
with $\epsilon_j=\epsilon(\psi_j)$, $(Tf,\psi_j)=b_j \theta_j$ and $\theta_j=(f,\phi_j)$.  This model is often considered in the inverse problem literature, see eg \cite{MR2421941}. Setting $Y_j=b_j^{-1} Z_j$ and $\sigma_j=\sigma b_j^{-1}$ for all $j\in \N^{\star}$, we obtain (\ref{mod1}). Hence inference on observations from model (\ref{mod1}) provides the same results for inverse problems. We stress that if estimation issues for inverse problem have been well studied over the past years (see for instance \cite{osulliv}, \cite{MR2421941} or  \cite{loublud1} for a model selection approach),  tests for inverse problems have either only been tackled, in the general case, by very few authors among them  we refer to \cite{munkpreprint},  \cite{IngSap}, \cite{autre} or have been investigated only  for the  very specific case of the convolution problem, see in \cite{Butucea}  and references therein. \vskip .1in

For all $\theta \in l_2(J)$, we set $\|\theta\|_2^2=\sum_{j \in J} \theta_j^2$ and for all $\theta\in l_{\infty}(J)$, $\| \theta\|_{\infty} = \sup_{j\in J} |\theta_j | $. The purpose of this paper is to  provide rates of testing for the hypothesis $"\theta=0"$ against the alternative $"\|\theta\|_q\geq \rho "$ where $q=2$ or $q=\infty$ . More precisely, given $\alpha, \beta \in ]0,1[$, a level $\alpha$ test $\Phi_{\alpha}$ of the null hypothesis $"\theta=0"$, and a class of vectors $\F \subset l_q(J)$, we define the uniform separation rate $\rho_q(\Phi_{\alpha},\F,\beta)$ of the test $\Phi_{\alpha}$
over the class $\F$ with respect to the $l_q$ norm  as the smallest radius $\rho$ such that 
the test guarantees a power greater that $1-\beta$
for all alternatives $\theta \in \F$ such that $\|\theta\|_q\geq \rho$. More formally
$$ \rho_q(\Phi_{\alpha},\F,\beta)=\inf \ac{\rho>0, \inf_{\theta \in \F, \|\theta\|_q \geq \rho}
\PP_{\theta}(\Phi_{\alpha} \mbox{ rejects } ) > 1-\beta} .$$
Let us now define $ \rho_q(\F,\alpha, \beta)$ as the infimum over all level $\alpha$ test $\Phi_{\alpha}$ of the quantity $\rho_q(\Phi_{\alpha},\F,\beta)$. This quantity will be called the $(\alpha,\beta)$ minimax rate of testing over the class $\F$. The aim of the paper is to determine this minimax rate of testing over various classes of alternatives $\F$, for the test of null hypothesis "$\theta =0$" in Model (\ref{mod1}) with respect to the $l_2$ and $l_{\infty}$ norms.\vskip .1in

The main reference for computing minimax rates of testing over non parametric alternatives is the series of paper due to Ingster\cite{Ing}, where various statistical models and a wide range of sets of alternatives are considered. Lepski and Spokoiny \cite{Lepspok} obtained minimax rates of testing over Besov bodies $\mathcal{B}_{s,p,q }(R)$ in the irregular case (when $0<p<2$),  see also Ingster and Suslina \cite{IngSuslina}. Ermakov \cite{Erm} determines a family of asymptotic minimax tests  for testing that the signal belongs to a parametric set  against  nonparametric sets of alternatives in the heteroscedastic Gaussian white noise.  In all these references, asymptotic minimax rates of testing are established. 
In Model (\ref{mod1}), with $\sigma_j=\sigma$ for all $j \in J$,  Baraud \cite{yannick} consider a non asymptotic point of view, which means that the noise level $\sigma$ is not assumed to converge towards $0$. This is the point of view that we adopt in this paper. We give a precise expression of the dependency of the minimax rates of testing with respect to the sequence $(\sigma_j)_{j\in J}$. The particular cases of interest correspond to polynomial and exponentially increasing sequences, which in the case of Model~(\ref{mod2}) leads to the so-called mildly and severely ill-posed inverse problems. We do not aim at providing the adaptive minimax rates of testing, which will be the core of a future work.\vskip .1in

The paper is organised as follows. In Section 2, we provide lower bounds for the minimax separation rate over classes of vectors $\theta$ with a finite number of non-zero coefficients, which yet covers sparse signals. In Section 3, we determine upper bounds for those minimax rates. In Section 4, we compute minimax rates of testing over ellipsoids and $l_p$ balls. The proofs are gathered in Sections 5 and 6. \\

To end this introduction, let us define some notations. Let  $ Y=(Y_j)_{j \in J}$ obey to Model (\ref{mod1}). We denote by $\theta$ the vector (or sequence) $(\theta_j)_{j \in J}$ and by $\PP_{\theta}$ the distribution of $Y$. All along the paper, we consider the test of null hypothesis "$\theta=0$".
Let $\alpha \in ]0,1[$ be some prescribed level. A test function $\Phi_{\alpha}$ is a measurable function of the observation $Y$, with values in $\ac{0,1}$. The null hypothesis is accepted if $\Phi_{\alpha}=0$ and rejected if $\Phi_{\alpha}=1$. Finally, for all $x\in \mathbb{R}$, we denote by $\lfloor x \rfloor$ the greater integer smaller than $x$ and we set $\lceil x \rceil = \lfloor x \rfloor +1$.

\section{Lower bounds} 
The bounds will be established for two classes of signals characterized by their non zeros coefficients. The first one deals with the elementary case where the coefficients are equal to zero after a certain rank. The second one concerns the so-called sparse signals which are defined by the amount of non zeroes coefficients which can be located at different scales.
\subsection{Lower bounds in ${l_2}$ norm}

In this section, we generalize the  results obtained by Baraud \cite{yannick} in an homoscedastic model to the heteroscedastic Model (\ref{mod1}).  \\
 
 We first give a lower bound for the minimax rate of testing over the set $S_D$, defined for all $D\geq 1$ by 
$$ S_D=\ac{ \theta \in l_2(J), \forall j>D, \theta_j=0 }.$$
When $J=\{1, \ldots,N\}$, we assume that $D \leq N$. 

\begin{prop} \label{borninf1}
Assume that $ Y=(Y_j)_{j \in J}$ obeys to Model (\ref{mod1}). 
Let $\beta \in ]0, 1-\alpha[$, $c(\alpha, \beta)= \pa{2\ln(1+4(1-\alpha-\beta)^2)}^{1/2} $ and 
$$\rho_D^2= c(\alpha,\beta) \pa{\sum_{j=1}^D \sigma_j^4}^{1/2}.$$
The following result holds:
$$ \forall \rho \leq \rho_D, \ \inf_{\Phi_{\alpha}} \sup_{\theta\in S_D, \|\theta\|_2 =\rho} \PP_{\theta}(\Phi_{\alpha}=0) \geq \beta.$$
This implies that the minimax rate of signal detection over $S_D$ with respect to the $l_2$ norm satisfies
$$ \rho_2(S_D, \alpha, \beta) \geq \rho_D.$$
\end{prop}
This proposition can be understood as follows: whatever the $\alpha$-level test chosen, for all $\rho \leq \rho_D$, there exists a signal $\theta \in S_D$ with 
norm $\rho$ such that the error of the second kind is greater than $\beta$. The results obtained in Proposition \ref{borninf1} coincide with the lower bound established by Baraud \cite{yannick} in the homoscedastic model ($\sigma_j=\sigma \quad \forall j\in J$).\\

Let us now consider the problem of sparse signal detection. Let $k, n \in \mathbb{N}^*$ with $k\leq n$. When $J= \ac{1,\ldots, N}$, we assume that $n \leq N$. We want to obtain lower bounds for the minimax separation rate of signal detection over the set $\cS_{k,n}$ defined by
\beq \label{Skn}
\cS_{k,n}=\ac{\theta \in l_2(J),\ \forall j>n,\  \theta_j=0, \mbox{ Card } \ac{j\leq n, \theta_j\neq 0}\leq k}. 
\eeq

\begin{theo}\label{minodure}
 Assume that $ Y=(Y_j)_{j \in J}$ obeys to Model (\ref{mod1}). Let $\sigma_{(1)}\leq \sigma_{(2 )}\leq \ldots \leq \sigma_{(n)}$, we define for all $l \in \{0, \ldots, n-k\}$,
 \beq \label{somvar} 
 \Sigma_{l,k}^2=\sum_{j=l+1}^{l+k} \sigma_{(j)}^2 . 
 \eeq
 Let $\beta \in ]0, 1-\alpha[$, such that $\alpha+\beta \leq 59 \%$. Let 
 \beq \label{rhosauvage}
\rho_{k,n}^2=\cro{ \max_{0\leq l \leq n-k}  \Sigma_{l,k}^2 \ln\pa{1+\frac{n-l}{k^2}\vee \sqrt{\frac{n-l}{k^2}}}\vee \pa{\sum_{j=n-k+1}^n  \sigma_{(j)}^4}^{1/2}}. 
\eeq
For all level $\alpha$ test $\Phi_{\alpha}$, 
$$ \inf_{\theta \in \cS_{k,n}, \|\theta\|_2 \geq \rho_{k,n}} \PP_{\theta} (\Phi_{\alpha}=1) \leq 1-\beta.$$
This implies that the minimax rate of signal detection over $ \cS_{k,n}$ with respect to the $l_2$ norm satisfies 
$$ \rho_2(\cS_{k,n}, \alpha, \beta) \geq \rho_{k,n}.$$
\end{theo}

{\bf Comments :} Let us consider three cases governing the behaviour of the sequence $(\sigma_j)_{j \in J}$.
\bee
\item In the homoscedastic case, $\sigma_j = \sigma$ for all $j \in J$. In this case, $\Sigma_{l,k}^2= \sigma^2 k$ for all $l$ and, taking $l=0$, we obtain that
$$ \rho_{k,n}^2 \geq \sigma^2 k\ln\pa{1+\frac{n}{k^2}\vee \sqrt{\frac{n}{k^2}}},$$
which corresponds to the lower bound established by Baraud \cite{yannick}. 
\item When $k\leq n/2$ and $ \Sigma_{\lfloor n/2 \rfloor,k}^2 \geq C \Sigma_{n-k,k}^2$ for some absolute constant $C$ (independent of $k$ and $n$), we obtain that
\beq \label{minopoly}
\rho_{k,n}^2 \geq   \cro{C\Sigma_{n-k,k}^2 \ln\pa{1+\frac{n}{2 k^2}\vee \sqrt{\frac{n}{2 k^2}}}
\vee \pa{\sum_{j=n-k+1}^n  \sigma_{(j)}^4}^{1/2}}.
\eeq
At the price of a factor $2$ in the logarithm ($n$ is replaced by $n/2$), the variance term appearing in the lower bound for $\rho_{k,n}^2$ is $ \Sigma_{n-k,k}^2$ which corresponds to the largest possible variance for a set of cardinality $k$ in $\{1, \ldots ,n\}$, indeed
$$ \Sigma_{n-k,k}^2 = \max_{m \in \M_{k,n}} \sum_{j \in m} \sigma_j^2,$$
where $ \M_{k,n}$ denotes the set of all subsets of  $\{1, \ldots ,n\}$ with cardinality $k$. \\
This situation occurs for example when $(\sigma_j)_{j\in J} $ grows at a polynomial rate : $ \sigma_j=j^{2\gamma}$ for some $\gamma >0$. In this case,
$$\Sigma_{n-k,k}^2 \leq k n^{2\gamma}, \quad  \Sigma_{n/2,k}^2 \geq k \pa{\frac{n}{2}}^{2\gamma}\geq \frac1{2^{2\gamma}} \Sigma_{n-k,k}^2 .  $$
\item When  $(\sigma_j)_{j\in J} $ grows at an exponential rate : $\sigma_{j}=\exp(\gamma j)$ for some $\gamma >0$, we obtain that $\rho_{k,n}^2 \geq \sigma_{(n)}^2$
\eee
\subsection{Lower bounds in ${l_{\infty}}$ norm}
In this section, we provide lower bounds for the minimax rate of signal detection in Model (\ref{mod1}) with respect to the $l_{\infty}$ norm. These bounds can be derived from the Theorem \ref{minodure}. Indeed, let  $\mathcal{S}_{k,n}$ be defined by (\ref{Skn}).
For $\theta \in \mathcal{S}_{1,n}$, $\|\theta\|_2=\|\theta\|_{\infty}$. Hence, denoting by $\rho_{\infty}(\mathcal{F},\alpha,\beta) $ the $(\alpha, \beta)$  minimax rate of testing over the class $ \mathcal{F}$ with respect to the $l_{\infty}$ norm, we obtain that
$$ \rho_{\infty}(\mathcal{S}_{1,n},\alpha,\beta) =\rho_2(\mathcal{S}_{1,n},\alpha,\beta) .$$
By definition of the $(\alpha, \beta)$  minimax rate of testing it is obvious that if $\mathcal{S}\subset \mathcal{S'}$, then $\rho_{\infty}(\mathcal{S})
\leq \rho_{\infty}(\mathcal{S}')$. $\mathcal{S}_{1,n} \subset \mathcal{S}_{k,n} $ for all $1\leq k\leq n$. This implies that
for all $1\leq k\leq n$,
$$ \rho_{\infty}(\mathcal{S}_{k,n},\alpha,\beta)\geq \rho_2(\mathcal{S}_{1,n},\alpha,\beta).$$
This leads to the following corollary of Theorem \ref{minodure} :
\begin{cor} \label{borneinfinfinie}
Assume that $ Y=(Y_j)_{j \in J}$ obeys to Model (\ref{mod1}). 
 Let $\beta \in ]0, 1-\alpha[$, such that $\alpha+\beta \leq 59 \%$. Let
  $$\cS_{k,n}=\ac{\theta \in l_2(J),\ \forall j>n,\  \theta_j=0, \mbox{ Card } \ac{j\leq n, \theta_j\neq 0}\leq k}. $$
 Let $\sigma_{(1)}\leq \sigma_{(2)}\leq \ldots \leq  \sigma_{(n)}$. 
 We define 
 \beq \label{rhosauvageinf}
\rho_{n,\infty}=\max_{0\leq l \leq n-1}  \sigma_{(l+1)}\sqrt{ \ln\pa{1+{n-l} }}.
\eeq
The following result holds :
$$\forall 1\leq k\leq n,  \quad \rho_{\infty}(\cS_{k,n}, \alpha, \beta) \geq \rho_{n,\infty}.$$
\end{cor}
The proof of the corollary follows directly from the arguments given above and therefore will be omitted. \\
{\bf Comments :} 
\bee
\item Let $ S_n=\ac{ \theta \in l_2(J), \forall j>n, \theta_j=0 }.$
Note that $  S_n = \cS_{n,n}$, hence 
$\forall n\geq 1, \quad \rho_{\infty}(S_n, \alpha, \beta) \geq \rho_{n,\infty}.$
\item When $\sigma_j=\sigma $ for all $1\leq j\leq n$, we obtain $ \rho_{n,\infty}=\sigma \sqrt{ \ln\pa{n+1}}.$
\item When $(\sigma_j)_{j\in J} $ grows at a polynomial rate : $ \sigma_j=j^{\gamma}$ for some $\gamma >0$, we obtain that  $ \rho_{n,\infty} \geq C(\gamma) n^{\gamma} \sqrt{\ln(n)}$ by taking $l=\lfloor n/2 \rfloor$ in (\ref{rhosauvageinf}). 
\item When $(\sigma_j)_{j\in J} $ grows at an exponential  rate : $ \sigma_j=\exp(j\gamma)$ for some $\gamma >0$, we obtain that $ \rho_{n,\infty} \geq  C \exp(n\gamma)$ for some constant $C$. 
\eee

\section{Upper bounds}
In this section, we give upper bounds for the $(\alpha,\beta)$ minimax rates  of testing over the sets $ S_D$ and $\cS_{k,n}$  that we compare with the lower bounds obtained in the previous section. \\
In order to show that the $(\alpha,\beta)$ minimax rate of testing with respect to the $l_q$ norm 
over a set $\F$ is bounded from above by $\rho$, it suffices to define a test statistic $\Phi_{\alpha}$ such that  the power of the test at each point  $\theta$ in $\F$ satisfying $\|\theta\|_q\geq \rho$ is greater than $1-\beta$. 

\begin{prop}\label{bornesup1}
Assume that $ Y=(Y_j)_{j \in J}$ obeys to Model (\ref{mod1}). 
Let $\alpha ,\beta  \in ]0, 1[$, and let $t_{D, 1-\alpha}(\sigma) $ denote the $1-\alpha $ quantile of $ \sum_{j=1}^D \sigma_j^2 \epsilon_j^2$ :
$$  \PP\pa{ \sum_{j=1}^D \sigma_j^2 \epsilon_j^2 \geq t_{D, 1-\alpha}(\sigma)} =\alpha.$$
Let $\Phi_{\alpha} $ be the test defined by
\begin{equation}
\Phi_{\alpha} =\un_{ \sum_{j=1}^D Y_j^2 > t_{D, 1-\alpha}(\sigma)} .
\label{eq:test1}
\end{equation}
Then, $\Phi_{\alpha} $ is a level-$\alpha$ test :
$$\PP_{\theta=0} (\Phi_{\alpha} =1)=\alpha. $$
Moreover, there  exists an absolute constant $C$  such that for all $ \theta \in S_D$, setting
$$C(\alpha, \beta)=C(\log(\beta^{-1}) + \log(\alpha^{-1})+\sqrt{\log (\alpha^{-1})} +\sqrt{\log (\beta^{-1})})$$
$$ \|\theta\|_2^2 \geq C(\alpha, \beta)  \pa{\sum_{j=1}^D \sigma_j^4}^{1/2} \Rightarrow \PP_{\theta} (\Phi_{\alpha} =1) > 1-\beta .$$
Hence, we obtain that
$$ \rho_2^2(S_D, \alpha, \beta) \leq  C(\alpha, \beta)  \pa{\sum_{j=1}^D \sigma_j^4}^{1/2}. $$
\end{prop}
Note that this bound coincides to the upper bound found in  Proposition~\ref{borninf1}. Hence, it proves that this lower bound is sharp. \\
\ \\
Let us now propose a testing procedure for sparse signal detection. This procedure will be defined by a combination of two tests. The first one is based on a thresholding method, which was already used for detection of irregular alternatives in Baraud et al \cite{BHL} and in Fromont et al \cite{FLR}. The second one is the test considered in Proposition \ref{bornesup1} with $D=n$, which will be powerful when  $k$ is larger that $\sqrt{n}$. 
\begin{theo}\label{bornesup2}

Assume that $ Y=(Y_j)_{j \in J}$ obeys to Model (\ref{mod1}). 
Let $\alpha ,\beta  \in ]0, 1[$, and let $t_{n, 1-\alpha}(\sigma) $ denote the $1-\alpha $ quantile of $ \sum_{j=1}^n \sigma_j^2 \epsilon_j^2$.
Let $\Phi_{\alpha}^{(1)} $ be the test defined by
$$  \Phi_{\alpha}^{(1)} =\un_{ \sum_{j=1}^n Y_j^2 > t_{n, 1-\alpha}(\sigma)} .$$
Let $q_{n, 1-\alpha} $ denote the $1-\alpha $ quantile of $ \max_{1\leq j\leq n } \epsilon_j^2$.
Let $\Phi_{\alpha}^{(2)} $ be the test defined by
\beqe
\Phi_{\alpha}^{(2)} &=& 1 \mbox{ if } \max_{1\leq j\leq n} \pa{ \frac{Y_j^2}{\sigma_j^2} } > q_{n, 1-\alpha}\\
&=&0 \mbox{ otherwise. }
\eeqe
We define $ \Phi_{\alpha} =\max\pa{ \Phi_{\alpha/2}^{(1)}, \Phi_{\alpha/2}^{(2)}}$.
Then, $\Phi_{\alpha} $ is a level-$\alpha$ test :
$$\PP_{\theta=0} (\Phi_{\alpha} =1)=\alpha. $$
There  exists a constant $C(\alpha, \beta) $ such that for all $ \theta \in \mathcal{S}_{k,n}$ satisfying
\beq \label{majindividuelle}
 \|\theta\|_2^2 \geq C(\alpha, \beta)  \cro{ \pa{\sum_{j=1}^n \sigma_j^4}^{1/2} \wedge \pa{\sum_{j,\theta_j\neq 0} \sigma_j^2} \ln(n)},
 \eeq
we have 
$$ \PP_{\theta} (\Phi_{\alpha} =1) > 1-\beta .$$
Hence, we obtain that
\beq \label{majosauvage}
\rho_2^2(\mathcal{S}_{k,n},\alpha, \beta) \leq  C(\alpha, \beta) \cro{ \pa{\sum_{j=1}^n \sigma_j^4}^{1/2} \wedge \Sigma_{n-k,k}^2 \ln(n)}, 
\eeq
where $\Sigma_{l,k}^2$ has been defined in (\ref{somvar}). 
\end{theo}
{\bf Comments :} Let us compare these results with the lower bounds obtained in Theorem \ref{minodure}. 
\bee 
\item We first assume that  $(\sigma_j)_{j\in J} $ grows at a polynomial rate :
$\forall j \in J$, $ \sigma_j =  \sigma j^{\gamma}$ for some $\gamma \geq 0$ (this includes the homoscedastic case). In this case,  when $k\leq n/2$ there exists a constant $C>0$ such that $ \Sigma_{\lfloor n/2 \rfloor,k}^2 \geq C\Sigma_{n-k,k}^2$. A lower bound for the $(\alpha,\beta)$ minimax separation rate of signal detection over $\mathcal{S}_{k,n}$ is given by (\ref{minopoly}). This lower bound has to be compared with the upper bound (\ref{majosauvage}).\\
\bei
\item When $k=n^l$ with $l<1/2$, the upper and lower bounds coincide and are of order $ \Sigma_{n-k,k}^2 \ln(n)$. 
\item  When $k=n^l$ with $l\geq 1/2$, the lower bound is of order $ \Sigma_{n-k,k}^2  {\sqrt{n}}/{k}$ and 
$\Sigma_{n-k,k}^2 \geq C k \sigma^2 n^{2\gamma}$, which leads to a lower bound of order $ C \sigma^2 n^{2\gamma+1/2}$. The upper bound is smaller that $ \pa{\sum_{j=1}^n \sigma_j^4}^{1/2} $, which is smaller than $ \sigma^2 n^{2\gamma+1/2}$. Hence, the two bounds coincide.
\item When $ k=\sqrt{n}/\phi(n)$ where $\phi(n) \rightarrow +\infty$  and $  \phi(n) /n \rightarrow 0 $ as $n\rightarrow +\infty$ (typically $\phi(n)=\ln(n)$), the lower bound is of order $ \Sigma_{n-k,k}^2 \ln(\phi(n))$ and the upper bound  is of order $ \Sigma_{n-k,k}^2 \ln(n)$. In this case, the upper and lower bound do not coincide, up to a logarithmic term. 
\eei
\item Let us now assume that  $(\sigma_j)_{j\in J} $ grows at an exponential rate : $\forall j \in J$, $ \sigma_j =\sigma \exp(\gamma j) $ for some $\gamma > 0$. The lower bound is greater than $\sigma_n^2= \sigma^2 \exp(2 \gamma n)$ and the upper bound is smaller that $  C(\alpha, \beta)\pa{\sum_{j=1}^n \sigma_j^4}^{1/2}$, which is bounded from above by $ C(\alpha, \beta,\gamma)\sigma^2 \exp(2 \gamma n)$. Hence the two bounds coincide. Note that in this case, the test  $ \Phi_{\alpha}^{(2)} $ based on thresholding  is useless and one can simply consider that test 
$$ \Phi_{\alpha} =\Phi_{\alpha}^{(1)},$$
which achieves the lower bound for the separation rate. 
\item The result stated in (\ref{majindividuelle}) is more precise than the minimax upper bound given in (\ref{majosauvage}). If the set $ J_1= \ac{j, \theta_j\neq 0 }$ corresponds to small values for the variances $(\sigma_j)_{j\in J_1}$, it is not required that $\|\theta\|_2^2$ is greater than the right hand term in (\ref{majosauvage}) for the test to be powerful for this value of $\theta$. The minimax bound given in (\ref{majosauvage}) corresponds to the worst situation, that is the case where the set $ J_1$ corresponds to the largest values for the variances. 
\eee
Let us now present an upper bound for the minimax separation rate with respect to the $l_{\infty}$ norm. 

\begin{cor} \label{bornesupinfinie}

Assume that $ Y=(Y_j)_{j \in J}$ obeys to Model (\ref{mod1}). Let $\sigma_{(1)} \leq \ldots \leq \sigma_{(n)}$. 
Let $\alpha ,\beta  \in ]0, 1[$. 
Let $\Phi_{\alpha} $ be the test defined in Theorem \ref{bornesup2}. 
There exists a constant $C(\alpha,\beta)$ such that for all $k \in \{1,\ldots, n\}$, for all $\theta\in \mathcal{S}_{k,n}$ such that
$$\|\theta\|_{\infty} \geq C(\alpha,\beta) \pa{\sigma_{(n)} \sqrt{\ln(n)}  \wedge \pa{ \sum_{j=1}^n \sigma_j^4}^{1/4}}$$
we have 
$$\PP_{\theta} (\Phi_{\alpha} =1) > 1-\beta .$$
This implies that for all $k \in \{1,\ldots, n\}$,
$$ \rho_{\infty}(\mathcal{S}_{k,n},\alpha,\beta)\leq C(\alpha,\beta) \pa{\sigma_{(n)} \sqrt{\ln(n)}  \wedge \pa{ \sum_{j=1}^n \sigma_j^4}^{1/4}}  $$
\end{cor}

This upper bound coincides with the lower bound obtained in Corollary \ref{borneinfinfinie} when the sequence $(\sigma_j)_{j\in J} $ is constant or grows at a polynomial or at an exponential rate.

\section{Minimax rates over ellipsoids and $l_p$ balls}
In the previous sections, the only constraint on the signal was expressed through the number of non-zero coefficients. In several situations, one deals instead with infinite sequences having a finite number of significant coefficients, the reminder being considered as negligible (in a sense which will be precised later on). To this end, we consider in this section a slightly different framework. Our aim is to study the link between the decay of the $\theta_k$'s and the associated rate of testing. We consider in the following two different kinds of function spaces: ellipsoids and $l_p$-bodies.

\subsection{Minimax rates of testing over ellipsoids}
In the following, we assume that the sequence $\theta=(\theta_j)_{j\in J}$ belongs to the ellipsoid $\mathcal{E}_{a,2}(R)$ defined as
$$ \mathcal{E}_{a,2}(R) = \left\lbrace \nu \in l_2(J),  \sum_{j\in J} a_j^2 \nu_j^2 \leq R^2 \right\rbrace, $$
where $a=(a_k)_{k\in J}$ denotes a monotone non-decreasing sequence. For instance, if $\theta$ corresponds to the sequence of Fourier coefficients of a function $f$ and $a_j$ is of order $j^s$ with $s>0$, then assuming that $\theta\in \mathcal{E}_{a,2}(R)$ is equivalent to impose conditions on the $s$-th derivative of $f$. The belonging to $\mathcal{E}_{a,2}(R)$ may be seen as a regularity assumption on our signal. The following result characterises the minimax rate of testing over $\mathcal{E}_{a,2}(R)$.

\begin{prop} 
\label{thm:minimax_ell}
Let $\alpha,\beta$ be fixed and denote by $\rho_2(\mathcal{E}_{a,2}(R),\alpha,\beta)$ the minimax rate of testing over $\mathcal{E}_{a,2}(R)$ with respect to the $l_2$ norm. Then
$$ \rho_2^2(\mathcal{E}_{a,2}(R),\alpha,\beta) \geq \sup_{D\in J} (\rho_D^2 \wedge R^2 a_D^{-2}),$$
where $\rho_D^2$ has been introduced in Proposition \ref{borninf1}. Moreover, for all $D\in J$,
$$ \sup_{\theta \in \mathcal{E}_{a,2}(R), \| \theta \|_2^2 \geq C\rho_D^2 + R^2 a_D^{-2}} P_{\theta}(\Phi_{\alpha}^{}=0) \leq \beta,$$
where $C=C(\alpha,\beta)$ is a positive constant depending only on $\alpha$ and $\beta$ and $\Phi_{\alpha}$ denotes the test introduced in Proposition \ref{bornesup1}. Hence,
$$ \rho_2^2(\mathcal{E}_{a,2}(R),\alpha,\beta) \leq \inf_{D\in J} (C\rho_D^2 + R^2 a_D^{-2}),$$
\end{prop}
Proposition \ref{thm:minimax_ell} presents both an upper and a lower bound for the minimax rate of testing over $\mathcal{E}_{a,2}(R)$. Remark that the upper bound is attained by the test $\Phi_{\alpha}$ introduced in Proposition \ref{bornesup1} where only signals with a finite number of non-zero coefficients were considered. We do not use the whole sequence $(Y_j)_{j\in J}$ in order to test the null hypothesis $"\theta = 0 "$ but only the first $D$ coefficients. The price to pay is to introduce some bias in the testing procedure. However this bias can be controlled by taking advantage of the constraint expressed on the decay of $\theta$.

A good characterization of $\rho_2(\mathcal{E}_{a,2}(R),\alpha,\beta)$ can be obtained as soon as the lower and upper bounds are of the same order. As many statistical problems encountered in the literature, one has to find a trade-off between the bias $R^2 a_D^{-2}$ and some kind of variance term $\rho_D^2$. This trade-off can be performed in several situations, hence leading to explicit rates of convergence: see Corollary \ref{cor:1} below.\\
\\
Proposition \ref{thm:minimax_ell} presents the minimax rate of testing in a general setting. Several explicit rates can be obtained when introducing specific constraints on the sequences $(a_k)_{k\in J}$ and $(b_k)_{k\in J}$. These rates are summarized in Corollary \ref{cor:1}. Let $(\nu_k)_{k\in\N^{\star}}$ be a sequence real numbers. In the following, we write $\nu_k \sim k^{l}$ if there exist positive constants $c_1$ and $c_2$ such that, for all $k\in\N^{\star}$, $c_1 k^l \leq \nu_k \leq c_2 k^l$.

\newtheorem{corollary}{Corollary}
\begin{cor}
\label{cor:1}
Let $\alpha,\beta$ be fixed. We assume that $J=\mathbb{N}^{\star}$ and $(Z_j)_{j\in J}$ obeys to Model (\ref{mod2}). The table below presents the minimax rates of testing over the ellipsoids $\mathcal{E}_{a,2}(R)$ with respect to the $l_2$ norm. We consider various behaviours for the sequences $(a_k)_{k\in \N^{\star}}$ and $(b_k)_{k\in \N^{\star}}$. For each case, we give $f(\sigma)$ such that for all $1>\sigma>0$, $C_1(\alpha,\beta) f(\sigma) \leq \rho^2(\mathcal{E}_{a,2}(R),\alpha,\beta) \leq  C_2(\alpha,\beta) f(\sigma)$ where $C_1(\alpha,\beta)$ and $C_2(\alpha,\beta)$ denote positive constants independent of $\sigma$. \\

\begin{center}
\begin{tabular}{lcc}
\hline
                    & \textbf{Mildly ill-posed}               & \textbf{Severely ill-posed} \\
                    &    $b_k \sim k^{-t}$                    &    $b_k \sim \exp(-\gamma k^r)$ \\
\hline
$a_k \sim k^{s}$   &  $\sigma^{\frac{4s}{2s+2t+1/2}}$   & $ \left( \log(\sigma^{-2}) \right)^{-2s/r} $  \\
\hline
$a_k \sim \exp(\nu k^s)$ &  $\sigma^2 \left( \log(\sigma^{-2}) \right)^{(2t+1/2)/s}$  & $e^{-2\nu \tilde D^s}$ $(s\leq 1)$ \\
\hline
\end{tabular}
\end{center}
Here $\tilde D$ denotes the integer part of the solution of $\rho_D^2 = R^2 a_D^{-2}$.
\end{cor}

These rates have already been presented in the literature. The case $a_k \sim k^{s}$ and $b_k \sim k^{-t}$ was first studied in \cite{IngBook}. More recently, \cite{IngSap} deals with other cases. Similar rates are also available in \cite{Butucea} in the context of density estimation with errors in the variables. The aim of Corollary \ref{cor:1} is to show that our approach can lead to important minimax results. Our point of view in this paper is indeed entirely non-asymptotic and is not restricted to ill-posed inverse problems.

Concerning severely ill-posed problems with supersmooth functions (i.e. $b_k \sim \exp(-\gamma k^r)$ and $a_k \sim \exp(\nu k^s)$), we do not handle the general case since we assume that $s\leq 1$. When this assumption is violated, the upper and lower bounds in Proposition \ref{thm:minimax_ell} do not coincide: our test does not attain the minimax rate of testing. This is certainly due to our approach, which in some sense is related to a rough regularization scheme. We mention \cite{IngSap} for a complete study of this case.

\subsection{Minimax rates of testing over $l_p$-bodies with $0<p<2$}
\label{s:lpbodies}
Ellipsoids contain essentially smooth functions. In the particular case where $\theta$ corresponds to the Fourier coefficients of a given function $f$, the constraints expressed through the belonging to one of the spaces introduced above may be incompatible with the presence of discontinuities. In order to extend the covered cases, we consider in this subsection sequences $\theta$ belonging to $l_p$-bodies $\mathcal{E}_{a,p}(R)$ defined as
$$ \mathcal{E}_{a,p}(R) = \left\lbrace \nu \in l_2(J), \sum_{j\in J} a_j^p \nu_j^p \leq R^p \right\rbrace, $$
where $a=(a_k)_{k\in J}$ denotes a monotone non-decreasing sequence and $0<p<2$. The following theorem proposes a lower bound for the minimax rate of testing over such spaces. 

\begin{theo} 
\label{thm:lpbod_low}
Let $(Y_j)_{j\in J}$ obey to  Model~(\ref{mod1}). Let $\alpha,\beta$ be fixed and denote by $\rho_2(\mathcal{E}_{a,p}(R),\alpha,\beta)$ the minimax rate of testing over $\mathcal{E}_{a,p}(R)$ with respect to the $l_2$ norm.
For all $D\in J$ and for all $0 \leq l \leq D- \lceil \sqrt{D} \rceil$, we set
$$ \rho^2_{\lceil \sqrt{D} \rceil,D,l} = \Sigma_{l,\lceil \sqrt{D} \rceil}^2 \ln\pa{1+ \sqrt{1-\frac{l}{D}}},$$
where $\Sigma_{l,\lceil \sqrt{D} \rceil}^2$ is given in (\ref{somvar}). Then
$$ \rho^2(\mathcal{E}_{a,p}(R),\alpha,\beta) \geq \sup_{D\in J} \pa{\rho_1(D) \vee \rho_2(D)},$$
where
$$ \rho_1(D) = \max_{0\leq l \leq D- \lceil \sqrt{D} \rceil} \pa{ \sqrt{D}^{1-2/p} R^2 a_D^{-2} \frac{\Sigma_{l,\lceil \sqrt{D} \rceil}^2}{\Sigma_{D-\lceil \sqrt{D} \rceil,\lceil \sqrt{D} \rceil}^2} \wedge \rho^2_{ \lceil \sqrt{D} \rceil,D,l} },$$ and 
$$ \rho_2(D) = \sqrt{D}^{1-2/p} R^2 a_D^{-2} \wedge \left( \sum_{j=D-\lceil \sqrt{D} \rceil+1}^D \sigma_{(j)}^4 \right)^{1/2}.$$ 
\end{theo} 

To the end of this section, we assume that the sequence $(b_j)_{j\in J}$ is polynomially or exponentially increasing, which yet correspond to the main case of interest in inverse problems. The lower bounds in this particular setting are easier to handle, as proved in the following corollary.

\begin{cor}\label{corrr}
Let $(Z_j)_{j\in \mathbb{N^{\star}}}$ obey to the Model (\ref{mod2}) with $b_j=\sigma \sigma_j^{-1}$ for all $j\in \mathbb{N}^{\star}$. Then, assuming that for all $j\in \mathbb{N}^{\star}$, $\sigma_j=j^{\gamma}$ or that for all $j\in\mathbb{N}^{\star}$, $\sigma_j=\exp(\gamma j)$ for some $\gamma \geq 0$, we obtain
$$ \rho^2(\mathcal{E}_{a,p}(R),\alpha,\beta) \geq C(\gamma)\sup_{D\in J} (\rho_{\lceil \sqrt{D} \rceil, D}^2 \wedge \sqrt{D}^{1-2/p} R^2 a_D^{-2}) := \rho^2_{a,p,R},$$
where $\rho_{\lceil \sqrt{D} \rceil, D}^2$ is defined in (\ref{rhosauvage}).
\end{cor}
\ \\
In order to attain the lower bound presented above, a test similar to the one introduced in Proposition \ref{bornesup1} is not sufficient. On $l_p$-bodies, the bias after a given rank $D$ is indeed more difficult to control than for ellipsoids. Some significant coefficients (in a sense which will be precised in the proof) may be contained in the sequence $\theta$ after the rank $D$. Hence, we have to introduce specific tests in order to detect these coefficients. 

More precisely, for all $j\in J$ and $\alpha \in (0,1)$, introduce
$$ \Phi_{\lbrace j \rbrace,\alpha} = \mathbf{1}_{\left\lbrace |Y_j| \geq q_{j,\alpha} \right\rbrace},$$
where $q_{j,\alpha}$ denotes the $1-\alpha$ quantile of a Gaussian random variable with mean $0$ and variance $\sigma_j^2$. Then define
$$ \Phi_{\alpha}^{\dagger} = \Phi_{\mathrm{loc},\alpha/2} \wedge \Phi_{D^{\dagger},\alpha/2} \ \mathrm{with} \ \Phi_{\mathrm{loc},\alpha/2} = \sup_{j\in \lbrace D^{\dagger}+1,.. N \rbrace} \Phi_{\lbrace j \rbrace, 3\alpha/\pi^2(j-D^{\dagger})^2}, $$
where $\Phi_{D^{\dagger},\alpha/2}$ denotes the test constructed in Proposition \ref{bornesup1} and
\begin{equation}
D^{\dagger} = \inf \left\lbrace D\in J, R^2 a_D^{-2} \sqrt{D}^{1-2/p} \leq \rho_{\lceil \sqrt{D} \rceil, D}^2 \right\rbrace.
\label{Dagger}
\end{equation}
By convention, $D^{\dagger}=N$ if $J=\lbrace 1,\dots, N \rbrace$ and the set in (\ref{Dagger}) is empty. The following theorem emphasizes the performances of the test $\Phi_{\alpha}^{\dagger}$.

\begin{prop}
\label{thm:lpbod_upp}
Let $\alpha,\beta$ be fixed. We assume that the sequence $\left(a_j^{-p}b_j^{-(2-p)}\right)_{j\in\N^*}$ is monotone non-increasing.  Suppose that $J= \lbrace 1,\dots, N \rbrace$. Then
$$ \sup_{\theta \in \mathcal{E}_{a,p}(R), \| \theta \|^2 \geq \lambda_{\sigma} \rho^2_{a,p,R}} P_{\theta}(\Phi_{\alpha}^{\dagger}=0) \leq \beta,$$
with 
\begin{itemize}
\item $\lambda_{\sigma} = C \log(N)$ for mildly ill-posed problems, i.e $(\sigma_j)_{j\in J} \sim (j^t)_{j\in J}$ for $t>0$,
\item $\lambda_{\sigma} = C \log(N) \sqrt{D^{\dagger}}^{1-p/2}$ for severely ill-posed inverse problems, i.e. $(\sigma_j)_{j\in J} \sim (e^{\gamma j})_{j\in J}$ for $\gamma>0$,
\end{itemize}
where $C$ denotes a positive constant independent of $\sigma$. 
\end{prop}

Our test reaches the lower bound established in Proposition~\ref{thm:lpbod_low} up to a log term,  it is not sharp. Nevertheless this drawback is not characteristic of the heteroscedastic model  since a similar property occurs in the homoscedastic case: see \cite{yannick} for more details. Hence, the lower bound established in Proposition \ref{thm:lpbod_low} corresponds certainly to the minimax rate on $\mathcal{E}_{a,p}(R)$. 

For the sake of convenience, the upper bound is presented for $J= \lbrace 1,\dots, N \rbrace$ which, roughly speaking,  corresponds to the regression setting. Nevertheless, our result can be easily extended to the case where $J=\N^{\star}$. In such a situation, our test will be performed on $\lbrace 1,\dots, \tilde N \rbrace$, where $\tilde N$ is a trade-off between the bias after the rank $\tilde N$ on $\mathcal{E}_{a,p}(R)$ and the growth of $\log (N)$. A good candidate for $\tilde N$ is a power of $\sigma^{-2}$.

In order to conclude this discussion, we point out that we impose a condition on the sequence $\left(a_j^{-p}b_j^{-(2-p)}\right)_{j\in\N}$. This condition is necessary in order to control the bias after the rank $D^{\dagger}$. It always hold when $p=2$ since $(a_j)_{j\in\N}$ is an increasing sequence. When $p<2$, the considered function has to be sufficiently smooth with respect to the ill-posedness of the problem. A similar condition can be found for instance in \cite{Donoho}. \\

The corollary below deals with the particular case of mildly ill-posed problems with polynomial $l_p$-bodies, where the situation is easier to handle. 
\begin{cor}
Assume that $a_k \sim k^s$ and $b_k \sim k^{-t}$ for all $k\in \N^*$ where $s,t$ denote positive constants such that $s>t(2/p-1)$. Then
$$  C_2 \log(N) \sigma^{\frac{4s+2/p-1}{2s+2t+1/p}} \geq \rho^2(\mathcal{E}_{a,p}(R),\alpha,\beta) \geq C_1 \sigma^{\frac{4s+2/p-1}{2s+2t+1/p}},$$
where $C_1,C_2$ denote positive constant independent of $\sigma$.
\end{cor}
Remark that the sequence $\left(a_j^{-p}b_j^{-(2-p)}\right)_{j\in\N^*}$ is monotone non-increasing as soon as $s>t(2/p-1)$. Hence the conditions of Proposition \ref{thm:lpbod_upp} are satisfied. The proof follows the same argument as in Corollary \ref{cor:1}.

\section{Proofs}
\subsection{Proof of the lower bounds}
The proofs of the lower bounds use a Bayesian approach extending the methods developed in the papers by Ingster \cite{Ing} and by Baraud \cite{yannick}. We use the following lemma :
\begin{lemma}\label{lemgene}
Let $\F$ be some subset of $l_2(J)$. Let $\mu_{\rho}$ be some probability measure on 
$$\F_{\rho}= \ac{\theta  \in \F, \|\theta\|_2\geq \rho}$$
and let
$$ \PP_{\mu_{\rho}}=\int \PP_{\theta}d \mu_{\rho}(\theta).$$
Assuming that $\PP_{\mu_{\rho}}$ is absolutely continuous with respect to $\PP_0$, we define 
$$L_{\mu_{\rho}}(y)=\frac{d\PP_{\mu_{\rho}}}{d\PP_0}(y).$$
For all $\alpha >0$, $\beta\in ]0, 1-\alpha[$, if 
$$ \EE_0\pa{L^2_{\mu_{\rho^*}}(Y)} \leq 1+ 4(1-\alpha-\beta)^2,$$
then
$$ \forall \rho \leq \rho^*, \quad \inf_{\Phi_{\alpha}}\sup_{\theta \in \F_{\rho}}
\PP_{\theta}(\Phi_{\alpha}=0) \geq \beta .$$
This implies that
$$ \rho(\F, \alpha, \beta )\geq  \rho^*.$$
\end{lemma}
For the proof of this lemma, we refer to Baraud \cite{yannick}, Section 7.1. \\ 

\subsubsection{Proof of Proposition \ref{borninf1}}
Let $\rho >0$, we set for  $1\leq j \leq D$, 
$$ \theta_j={\omega_j}{\sigma_j^2}\rho\pa{\sum_{j=1}^D \sigma_j^4}^{-1/2}$$
where $(\omega_j, 1\leq j \leq D)$ are i.i.d. Rademacher random variables : $\PP(\omega_j=1)=\PP(\omega_j=-1)=1/2$. Let $\mu_{\rho}$ be the distribution of  $(\theta_1, \ldots, \theta_D)$.  $\mu_{\rho}$ is a probability measure on 
$$ \ac{\theta  \in S_D, \|\theta\|_2=\rho}. $$
Let us now evaluate the likelihood ratio $L_{\mu_{\rho}}(Y)=\frac{d\PP_{\mu_{\rho}}}{d\PP_0}(Y).$
\begin{eqnarray*}
L_{\mu_{\rho}}(Y) &=& \EE_{\omega} \cro{ \exp\pa{-\frac1{2}\sum_{j=1}^D \frac1{\sigma_j^2} \pa{Y_j-\frac{\sigma_j^2 \omega_j \rho}{\sqrt{\sum_{j=1}^D \sigma_j^4}}}^2} \exp\pa{\frac1{2}\sum_{j=1}^D \frac{Y_j^2}{\sigma_j^2}}}\\
 & =& \exp \pa{-\frac{\rho^2}{2} \frac{\sum_{j=1}^D \sigma_j^2}{\sum_{j=1}^D \sigma_j^4}} \prod_{j=1}^D \cosh \pa{\frac{\rho Y_j}{\sqrt{\sum_{j=1}^D \sigma_j^4}}}.
 \end{eqnarray*}
Let $Z$ be some standard normal variable. For all $\lambda \in \mathbb{R}$, 
\beq\label{cosh2}
\EE(\cosh^2(\lambda Z))= \exp(\lambda^2)\cosh(\lambda^2).
\eeq
Hence, since $Y_j/\sigma_j \sim \mathcal{N}(0,1)$, 
$$ \EE_{0}\pa{L^2_{\mu_{\rho}}(Y)} = \prod_{j=1}^D \cosh \pa{\frac{\rho^2 \sigma_j^2}{\sum_{j=1}^D \sigma_j^4}}.$$
Since for all $x \in \mathbb{R}$, $\cosh(x)\leq \exp(x^2/2)$, we obtain
$$ \EE_{0}\pa{L^2_{\mu_{\rho}}(Y)}\leq \exp\pa{\frac{\rho^4}{2\sum_{j=1}^D \sigma_j^4}}.$$
 For $\rho=\rho_D$  we obtain :
 $$ \EE_{0}\pa{L^2_{\mu_{\rho}}(Y)} \leq 1+4(1-\alpha-\beta)^2,$$
 which implies that $\rho(S_D, \alpha, \beta )\geq 
\rho_D$ by Lemma \ref{lemgene}.
\subsubsection{Proof of Theorem \ref{minodure}}
Without loss of generality, we can assume that the sequence $(\sigma_j)_{j \in J}$ is non decreasing (if this is not the case, we can reorder the observations $Y_j$). 
We fix some $l \in\{0,1,\ldots, n-k\}$. Let $\mathcal{M}_{k,l,n}$ denote the set of all subsets of $\ac{l+1,\ldots,n}$ with cardinality $k$.
 Let $\hat{m}$ be a random set of $\ac{l+1,\ldots,n}$, which is uniformly distributed on $\mathcal{M}_{k,l,n}$. This means that for all $m \in \mathcal{M}_{k,l,n}$, $\PP(\hat{m}=m)=1/C_{n-l}^k$.
 Let $(\omega_j, 1\leq j \leq n)$ be i.i.d. Rademacher random variables, independent of $\hat{m}$.  
  Let us recall that $$ \Sigma_{l,k}^2 =\sum_{j =l+1 }^{l+k} \sigma_j^2.$$
 We set
 \beq \label{minoalter}
 \theta_j=( \rho\omega_j  \sigma_j/{\Sigma_{l,k}} )\un_{j \in \hat{m}}
 \eeq
 Note that $\theta=(\theta_j)_{j \in J} \in \cS_{k,n}$ and that, since $(\sigma_j)_{j \in J} $ is non decreasing, 
 $$ \|\theta\|_2^2= \rho ^2 \frac{\sum_{j \in \hat{m}} \sigma_j^2}{\Sigma_{l,k}^2}\geq \rho^2.$$
\begin{eqnarray*}
 L_{\mu_{\rho}}(Y)&=&\EE_{\hat{m},\omega}\cro{ \exp\pa{-\frac1{2}\sum_{j\in J } \frac1{\sigma_j^2} \pa{Y_j-\theta_j}^2} \exp\pa{-\frac1{2}\sum_{j\in J} \frac{Y_j^2}{\sigma_j^2}}} \\
 &=&\EE_{\hat{m},\omega}\cro{ \exp\pa{ \sum_{j \in \hat{m}}\frac{Y_j\theta_j}{\sigma_j^2}}\exp\pa{-\frac12 \sum_{j \in \hat{m}} \frac{\theta_j^2}{\sigma_j^2 }}}.\\
 &=& \EE_{\hat{m},\omega}\cro{  \exp\pa{ \sum_{j \in \hat{m}} \frac{Y_j \omega_j\rho}{\sigma_j{\Sigma_{l,k}}}} \exp\pa{ 
 -\frac{k \rho^2}{2 \Sigma_{l,k}^2  }}}.
 \end{eqnarray*}
 \begin{eqnarray*}
 L_{\mu_{\rho}}(Y)&=& \frac{1}{C_{n-l}^k} \sum_{m \in \M_{k,l,n}} \EE_{\omega}\cro{\exp\pa{\sum_{j \in m} \frac{Y_j \omega_j \rho}{\sigma_j \Sigma_{l,k}}}\exp\pa{-\frac{k \rho^2}{2 \Sigma_{l,k}^2}}}\\
 &=& \exp\pa{-\frac{k \rho^2}{2 \Sigma_{l,k}^2}} \frac{1}{C_{n-l}^k} \sum_{m \in \M_{k,l,n}} \prod_{j \in m} \cosh\pa{\frac{\rho Y_j}{\sigma_j \Sigma_{l,k}}}.
 \end{eqnarray*}
 We use (\ref{cosh2}) together with $\EE(\cosh(\lambda Z))= \exp(\lambda^2/2)$ for a standard Gaussian variable $Z$. Since $Y_j/\sigma_j$ is a standard normal variable, we obtain that
 \beqe
  \EE_0 \pa{ L^2_{\mu_{\rho}}(Y)}&= &\exp\pa{-\frac{k \rho^2}{ \Sigma_{l,k}^2}}\frac1{(C_{n-l}^k)^2} \sum_{m,m' \in \M_{k,l,n}} \prod_{j \in m \backslash m'} \exp\pa{\frac{ \rho^2}{2 \Sigma_{l,k}^2}}\\
  &\times &\prod_{j \in m ' \backslash m} \exp\pa{\frac{ \rho^2}{2 \Sigma_{l,k}^2}} \prod_{j \in m \cap m'} \exp\pa{\frac{ \rho^2}{ \Sigma_{l,k}^2}} \cosh\pa{\frac{\rho^2}{\Sigma_{l,k}^2}}.
  \eeqe
Since, for all $  m,m' \in \M_{k,l,n}$, 
$$| m \backslash m' | +  | m' \backslash m| + 2| m \cap m'| = |m|+|m'|= 2k,$$
we obtain 
 $$ \EE_0 \pa{ L^2_{\mu_{\rho}}(Y)}= \frac1{(C_{n-l}^k)^2} \sum_{m,m' \in \M_{k,l,n}}\pa{\cosh\pa{\frac{\rho^2}{\Sigma_{l,k}^2}}}^{|m\cap m'|}. $$
 The end of the proof is similar to the proof of Theorem 1 in Baraud \cite{yannick}, similar arguments are also given in Fromont et al. \cite{FLR}. Let us recall these arguments.
 $$\EE_0 \pa{ L^2_{\mu_{\rho}}(Y)}=\EE  \cro{ \exp\pa{|\hm\cap \hm'|\ln  \cosh\pa{\frac{\rho^2}{\Sigma_{l,k}^2}}}},$$
 where $\hm,\hm'$ are independent random subsets with uniform distribution on $\M_{k,l,n}$.  
 For fixed $\hm$, $|\hm\cap \hm'|$ is an hypergeometric variable with parameters $(n-l,k,k/(n-l))$. We know from Aldous \cite{Aldous} that there exists a binomial variable $B$ with parameters $(k, k/(n-l))$ and a $\sigma-$ algebra $\mathcal{B}$  such that $ \EE ( B/\mathcal{B}) = |\hm\cap \hm'|  $. By Jensen's inequality, 
 $$ \EE_0 \pa{ L^2_{\mu_{\rho}}(Y)} \leq  \EE \cro{ \exp\pa{B \ln \cosh\pa{\frac{\rho^2}{\Sigma_{l,k}^2}}}}.$$
 Since $B$ is a binomial variable  with parameters $(k, k/(n-l))$,
 $$   \EE \cro{ \exp\pa{B \ln \cosh\pa{\frac{\rho^2}{\Sigma_{l,k}^2}}}} =  \exp \cro{ k \ln \pa{1 +\frac{k}{n-l}\pa{ \cosh\pa{\frac{\rho^2}{\Sigma_{l,k}^2}}-1}}} . $$
 Let $c = 1+4(1-\alpha-\beta)^2$, and $A= \frac{n-l}{k^2}\ln(c)$.  Since the function $\cosh$ is increasing on $\mathbb{R}^+$, 
 we obtain that if 
 $$ \frac{\rho^2}{\Sigma_{l,k}^2} \leq  \ln \pa{1+A+\sqrt{ 2A + A^2}},$$
 then
 $$ \cosh\pa{\frac{\rho^2}{\Sigma_{l,k}^2}} -1 \leq \frac{1}{2} \pa{A +\sqrt{2A+A^2}-1} + \frac{1}{2} \pa{A +\sqrt{2A+A^2}+1}^{-1} = A.$$
 We finally obtain that
 $$  \EE_0 \pa{ L^2_{\mu_{\rho}}(Y)} \leq  \exp \cro{ k \ln \pa{1 + \frac{k}{n-l} A } }\leq c.$$
 By Lemma \ref{lemgene}, this implies that 
 \beqe
  \rho(\cS_{k,n}, \alpha, \beta) &\geq & \Sigma_{l,k}^2 \ln \pa{1+A+\sqrt{ 2A + A^2}}   \\
  &\geq & \Sigma_{l,k}^2 \ln \pa{1+ 2A\vee \sqrt{ 2A }}.
  \eeqe
 If $\alpha+\beta \leq 0.59$,  $ \ln(c) \geq 1/2$, which implies that
 $$  \rho(\cS_{k,n}, \alpha, \beta) \geq \Sigma_{l,k}^2 \ln\pa{1+\frac{n-l}{k^2}\vee \sqrt{\frac{n-l}{k^2}}}.  $$
 Since this result holds for all $l\in\{0, n-k\}$, we get
$$\rho(\cS_{k,n}, \alpha, \beta) \geq \max_{0\leq l\leq n-k} \Sigma_{l,k}^2  \ln\pa{1+\frac{n-l}{k^2}\vee \sqrt{\frac{n-l}{k^2}}}.$$
In order to prove that 
$$\rho(\cS_{k,n}, \alpha, \beta) \geq \pa{\sum_{j=n-k+1}^n  \sigma_j^4}^{1/2},$$
we define, as in the proof of Proposition \ref{borninf1},
\begin{center}
\begin{tabular}{cll} 
 $\theta_j$&$=\omega_j \sigma_j^2 \rho\pa{\sum_{j=n-k+1}^n \sigma_j^4}^{-1/2} \quad$ &$\forall j\in\{n-k+1,\ldots,n\}$, \\
  & $=0 \quad $&$\forall j\notin\{n-k+1,\ldots,n\},$
 \end{tabular}
 \end{center}
where $(\omega_j, n-k+1\leq j \leq n)$ are i.i.d. Rademacher random variables. 
Note that $(\theta_j)_{j\in J} \in \cS_{k,n}$ and that $\|\theta\|_2^2 =\rho^2$. 
We now conclude as in the proof of Proposition \ref{borninf1}, using that $C(\alpha, \beta)=\sqrt{2 \ln (c)} \geq 1$.

\subsection{Proof of the upper bounds}
\subsubsection{Proof of Proposition \ref{bornesup1}}
In order to prove Proposition \ref{bornesup1}, we have to show that for all $\theta \in S_D$ such that
$\|\theta\|_2^2 \geq C(\alpha, \beta) \pa{\sum_{j=1}^D \sigma_j^4}^{1/2}$,
\beq \label{ine1}
\PP_{\theta} \left( \sum_{j=1}^D Y_j^2 \leq  t_{D, 1-\alpha}(\sigma) \right) < \beta .
\eeq
We denote  by $t_{D ,\beta}(\theta,\sigma)$ the $\beta $ quantile of $ \sum_{j=1}^D Y_j^2$, when  $ Y=(Y_j)_{j \in J}$ obeys to Model (\ref{mod1}). In order to prove (\ref{ine1}), it suffices to show that 
$$ t_{D, 1-\alpha}(\sigma)  < t_{D, \beta}(\theta,\sigma).$$
To prove this inequality, we will first give an upper bound for $ t_{D, 1-\alpha}(\sigma)$ and then a lower bound for $t_{D ,\beta}(\theta,\sigma)$. \\
{\bf Upper bound for $ t_{D, 1-\alpha}(\sigma)$ :} \\
We use an exponential inequality for chi-square distributions due to Laurent and Massart \cite{BLPM} (see Lemma 1). It follows from this inequality that for all $x \geq 0$, 
$$ \PP\pa{ \sum_{j=1}^D \sigma_j^2 (\epsilon_j^2 -1) \geq 2 \sqrt{x} \pa{\sum_{j=1}^D \sigma_j^4  }^{1/2}+ 2 x \sup_{1\leq j\leq D} (\sigma_j^2)  } \leq \exp(-x).$$
Setting $x_{\alpha}= \ln (1/\alpha)$, we obtain that
$$ t_{D, 1-\alpha}(\sigma) \leq \sum_{j=1}^D  \sigma_j^2 +  2 \sqrt{x_{\alpha}} \pa{\sum_{j=1}^D \sigma_j^4  }^{1/2} + 2 x_{\alpha} \sup_{1\leq j\leq D} (\sigma_j^2) .$$
Since $ \sup_{1\leq j\leq D} \sigma_j^2 \leq \pa{\sum_{j=1}^D \sigma_j^4  }^{1/2},$
\beq \label{quant1}
t_{D, 1-\alpha}(\sigma) \leq \sum_{j=1}^D  \sigma_j^2 +C(\alpha)\pa{\sum_{j=1}^D \sigma_j^4  }^{1/2}. 
\eeq 
\ \\
{\bf Lower bound for $ t_{D, \beta}(\theta,\sigma)$ :} \\
We prove the following lemma, which generalizes the results obtained by Birg\'e \cite{Lucien} to the heteroscedastic framework :
\begin{lemma}\label{khideux}
Let 
$$ Y_j=\theta_j + \sigma_j \epsilon_j,\quad 1\leq j \leq D,$$
where $ \epsilon_1, \ldots \epsilon_D$ are i.i.d. Gaussian variables with mean $0$ and variance $1$. \\
We define $\hT=\sum_{j=1}^D Y_j^2$ and 
$$\Sigma = \sum_{j=1}^D \sigma_j^4 +2 \sum_{j=1}^D \sigma_j^2 \theta_j^2 .$$
 The following inequalities hold for all $x \geq 0$ :
\beq \label{majo}
\PP\pa{\hT-\EE(\hT)\geq 2 \sqrt{\Sigma x} + 2 \sup_{1\leq j\leq D}(  \sigma_j^2 ) x }\leq \exp(-x).
\eeq
\beq \label{mino}
\PP\pa{\hT-\EE(\hT)\leq -2 \sqrt{\Sigma x}} \leq \exp(-x).
\eeq
\end{lemma}
The proof of this lemma is given in the Appendix. \\
Inequality (\ref{mino}) provides a lower bound for  $ t_{D, \beta}(\theta,\sigma)$. Indeed, setting $x_{\beta}= \log(1/\beta)$, we obtain that
$$  \PP\pa{\hT -\EE(\hT)\leq - 2 \sqrt {\Sigma x_{\beta}}} \leq \beta.$$
Hence, $ t_{D, \beta}(\theta,\sigma) \geq   \sum_{j=1}^D(\theta_j^2+\sigma_j^2) - 2 \sqrt {\Sigma x_{\beta}}$. 
(\ref{ine1}) is satisfied if $ t_{D, 1-\alpha}(\sigma) < t_{D ,\beta}(\theta,\sigma),$ which holds 
as soon as
 \beq \label{Condition}
 \sum_{j=1}^D\theta_j^2 - 2 \sqrt {\Sigma x_{\beta}} > 2 \sqrt{x_{\alpha}} \sqrt{\sum_{j=1}^D \sigma_j^4  }+ 2 x_{\alpha} \sup_{1\leq j\leq D} (\sigma_j^2).
 \eeq
Let us note that 
\beqe
\sqrt{\Sigma}= \sqrt{\sum_{j=1}^D \sigma_j^4+ 2 \sigma_j^2 \theta_j^2} &\leq& \sqrt{\sum_{j=1}^D \sigma_j^4}+ \sqrt{2} \sqrt{\sum_{j=1}^D \sigma_j^2\theta_j^2}\\
&\leq& \sqrt{\sum_{j=1}^D \sigma_j^4}+ \sqrt{2}  \sup_{1\leq j\leq D} (\sigma_j) \sqrt{\sum_{j=1}^D \theta_j^2}
\eeqe
Hence, the following inequality implies (\ref{Condition}) :
$$ \sum_{j=1}^D\theta_j^2  -2\sqrt{2} \sup_{1\leq j\leq D} (\sigma_j) \sqrt{x_{\beta}} \sqrt{\sum_{j=1}^D\theta_j^2 } -2 \sqrt{\sum_{j=1}^D \sigma_j^4}(\sqrt{x_{\beta}}+\sqrt{x_{\alpha}} )
-2 \sup_{1\leq j\leq D} (\sigma_j^2 )x_{\alpha} > 0.$$
This inequality holds if
$$  \sum_{j=1}^D\theta_j^2  \geq C \cro{\sup_{1\leq j\leq D} (\sigma_j^2) (x_{\beta}+x_{\alpha}) +\sqrt{\sum_{j=1}^D \sigma_j^4}( \sqrt{ x_{\beta}} +   \sqrt{ x_{\alpha}})},$$
where $C$ is an absolute constant (which can be taken equal to $8$).
Hence, we have proved that 
$$  \rho(S_D, \alpha, \beta) \leq C(\alpha, \beta) \sqrt{\sum_{j=1}^D \sigma_j^4},$$
which concludes the proof of Proposition \ref{bornesup1}.

\subsubsection{Proof of Theorem \ref{bornesup2}}
The test $\Phi_{\alpha}$ is obviously of level $\alpha$ thanks to Bonferroni's inequality :
\beqe
\PP_{0} (\Phi_{\alpha} =1) &\leq & \PP_{0} (\Phi_{\alpha/2}^{(1)}=1) + \PP_{0} (\Phi_{\alpha/2}^{(2)}=1)\\
&\leq& \frac{\alpha}{2}+ \frac{\alpha}{2} \leq \alpha.
\eeqe
Let us now evaluate the power of the test.
$$ \PP_{\theta} (\Phi_{\alpha} =1) \geq \max \pa{ \PP_{\theta} ( \Phi_{\alpha/2}^{(1)} =1)+ \PP_{\theta} ( \Phi_{\alpha/2}^{(2)} =1)}.$$
It follows from Proposition \ref{bornesup1} that for all $\theta \in \mathcal{S}_{k,n}$ such that
$$\|\theta\|_2^2 \geq C(\alpha,\beta) \pa{\sum_{j=1}^n \sigma_j^4}^{1/2},$$
we have $ \PP_{\theta} (\Phi_{\alpha}^{(1)} =1) > 1-\beta$. It remains to evaluate the power of the test $ \Phi_{\alpha/2}^{(2)}$. \\
\beqe
\PP_{\theta} (\Phi_{\alpha}^{(2)}=0) &=& \PP_{\theta} (\forall j \in \{1,\ldots,n\}, \frac{Y_j^2}{\sigma_j^2}\leq q_{n,1-\alpha}) \\
&\leq&\inf_{1\leq j\leq n} \PP_{\theta} (\frac{Y_j^2}{\sigma_j^2}\leq q_{n,1-\alpha}).
\eeqe
\beqe 
\PP_{\theta} (\frac{Y_j^2}{\sigma_j^2}\leq q_{n,1-\alpha}) &=& \PP \pa{ |\theta_j+\sigma_j \epsilon_j|
\leq \sigma_j \sqrt{q_{n,1-\alpha}}}\\
&\leq&  \PP \pa{ |\theta_j|-\sigma_j |\epsilon_j|\leq \sigma_j \sqrt{q_{n,1-\alpha}}}\\
&\leq&  \PP \pa{\sigma_j |\epsilon_j|\geq |\theta_j| -\sigma_j \sqrt{q_{n,1-\alpha}}}.
\eeqe
Let $q_{\beta}$ denote the $1-\beta$ quantile of $|\epsilon_j|$. We obtain that if 
\beq \label{condpuis}
\exists j \in \{1,\ldots,n\},\quad |\theta_j| > \sigma_j ( q_{\beta} + \sqrt{q_{n,1-\alpha}}),
\eeq
then 
$$ \PP_{\theta} (\Phi_{\alpha}^{(2)}=0) \leq \beta.$$
Condition (\ref{condpuis}) is equivalent to
$$\exists m \in \M_{k,n}, \sum_{j \in m} \theta_j^2 > \sum_{j \in m} \sigma_j^2 ( q_{\beta} + \sqrt{q_{n,1-\alpha}})^2.$$
In particular, if 
$$ \|\theta\|_2^2 > \pa{ \sum_{j,\theta_j\neq 0 } \sigma_j^2 }( q_{\beta} + \sqrt{q_{n,1-\alpha}})^2,$$
then (\ref{condpuis}) holds. 
This implies that for all $\theta \in \mathcal{S}_{k,n}$ such that 
$$\|\theta\|_2^2 > \max_{m \in \M_{k,n}}\pa{ \sum_{j \in m} \sigma_j^2 }( q_{\beta} + \sqrt{q_{n,1-\alpha}})^2,$$
we have $ \PP_{\theta} (\Phi_{\alpha}^{(2)}=0) < \beta.$
It remains to give an upper bound for  $q_{n,1-\alpha} $. We use the inequality $\PP(|\epsilon_1|\geq x)\leq \exp(-x^2/2)$. This leads to
\beqe
 \PP( \max_{1\leq j\leq n} \epsilon_j^2 \geq 2 \ln(n/\alpha)) &\leq & n \PP( |\epsilon_1|\geq  \sqrt{2 \ln(n/\alpha)}) \\
 &\leq &\alpha.
 \eeqe
Hence, $q_{n,1-\alpha} \leq 2 \ln(n/\alpha) $, which concludes the proof of Theorem \ref{bornesup2}.

\subsubsection{Proof of Corollary \ref{bornesupinfinie}}
Since $\|\theta\|_2^2 \geq \|\theta\|_{\infty}^2  \geq C^2(\alpha,\beta)\pa{\sum_{j=1}^n \sigma_j^4}^{1/2} $, we obtain from Theorem \ref{bornesup2} that 
$$  \PP_{\theta} (\Phi_{\alpha/2}^{(1)}=0) < \beta.$$
We obtained in the proof of Theorem \ref{bornesup2} that if (\ref{condpuis}) holds, 
then  $ \PP_{\theta} (\Phi_{\alpha}^{(2)}=0) < \beta.$
Since $q_{n,1-\alpha} \leq 2 \ln(n/\alpha) $, we obtain that there exists a constant $C(\alpha,\beta)$ such that  if $\|\theta\|_{\infty} \geq C(\alpha,\beta) \sigma_{(n)} \sqrt{\ln(n)}$, then $ \PP_{\theta} (\Phi_{\alpha/2}^{(2)}=0) < \beta.$

\subsection{Proof of minimax rates on ellipsoids and $l_p$-bodies}

\subsubsection{Proof of Proposition \ref{thm:minimax_ell}}

 We first prove the lower bound. For all $D\in J$, introduce $r_D^2 = \rho_D^2 \wedge R^2 a_D^{-2}$. Let $D$ be fixed. Then for all $\theta\in S_D$ such that $\|\theta\|_2^2 = r_D^2$
$$ \sum_{j\in J} a_j^2 \theta_j^2 = \sum_{j=1}^D a_j^2 \theta_j^2 \leq a_D^2 \| \theta \|_2^2 \leq R^2.$$
Hence
$$ \left\lbrace \theta\in S_D , \|\theta \|_2^2 = r_D^2 \right\rbrace \subset \left\lbrace \theta \in \mathcal{E}_{a,2}(R), \|\theta \|_2^2 \geq r_D^2 \right\rbrace.$$
Since $r_D\leq \rho_D $, we get from Proposition \ref{borninf1}
\begin{equation}
\inf_{\Phi_{\alpha}} \sup_{\theta \in \mathcal{E}_{a,2}(R), \|\theta \|_2 \geq r_D} P_{\theta} ( \Phi_{\alpha} =0) \geq \inf_{\Phi_{\alpha}} \sup_{\theta \in S_D, \|\theta \|_2 = r_D} P_{\theta} ( \Phi_{\alpha} =0) \geq \beta,
\label{eq:low1}
\end{equation}
where the infinimum is taken over all possible level-$\alpha$ testing procedures. Since inequality (\ref{eq:low1}) holds for all $D \in J$, we obtain $\rho^2(\mathcal{E}_{a,2}(R),\alpha,\beta) \geq \sup_{D \in J} (\rho^2_D \wedge R^2 a_D^{-2} )$. Concerning the upper bound, we know from Proposition \ref{bornesup1} that the test $\Phi_{\alpha}$ is powerful as soon as:
$$ \sum_{j=1}^{D} \theta_j^2 \geq C(\alpha,\beta) \rho_{D}^2 \Leftrightarrow \|\theta \|_2^2 \geq C(\alpha,\beta) \rho_{D}^2 + \sum_{k>D} \theta_k^2,$$
where $C(\alpha,\beta)$ denotes a positive constant. Since $\theta\in \mathcal{E}_{a,2}(R)$, we get
$$ \sum_{k>D} \theta_k^2 \leq R^2 a_{D}^{-2} \ \mathrm{and} \ \sup_{\theta\in \mathcal{E}_{a,2}(R),\| \theta \|^2 \geq C \rho_{D}^2 +R^2 a_{D}^{-2}} P_{\theta} (\Phi_{\alpha}^{} = 0) < \beta, $$
where $C=C(\alpha,\beta)$. This concludes the proof since the previous result holds for all $D\in J$.

\subsubsection{Proof of Corollary \ref{cor:1}} 
\textit{First case}: $a_k \sim k^s$ and $b_k \sim k^{-t}$. Choosing
$$\bar D = \left\lfloor \sigma^{\frac{2}{4s +4t +1}} \right\rfloor,$$
we can remark that $\rho_{\bar D}^2$ and $R^2 a_{\bar D}^{-2}$ are of the same order, hence leading to the desired rate.\\
\\
\textit{Second case}: $a_k \sim e^{\nu k^s}$ and $b_k\sim k^{-t}$. Set
$$ D_0 = \left\lceil \left( \frac{1}{2\nu} \log (\sigma^{-2}) \right)^{1/s} \right\rceil.$$
Then
\begin{eqnarray*}
\rho_2^2(\mathcal{E}_{a,2}(R),\alpha,\beta)  
& \leq & C\rho_{D_0}^2 + R^2 a_{D_0}^{-2},\\
& \leq & C \sigma^2 \left(\log(\sigma^{-2}) \right)^{(2t+1/2)/s} + \sigma^2 \leq (C+1) \sigma^2 \left(\log(\sigma^{-2}) \right)^{(2t+1/2)/s},
\end{eqnarray*}
where $C$ denotes a constant independent of $\sigma$. Concerning the lower bound, we set
$$ D_1 = \left\lfloor \left( \frac{1}{4\nu} \log (\sigma^{-2}) \right)^{1/s} \right\rfloor.$$
Then
\begin{eqnarray*}
\rho_2^2(\mathcal{E}_{a,2}(R),\alpha,\beta)
& \geq & \rho_{D_1}^2 \wedge R^2 a_{D_1}^{-2},\\
& \geq & C \sigma^2 \left(\log(\sigma^{-2}) \right)^{(2t+1/2)/s} \wedge \sigma = C \sigma^2 \left(\log(\sigma^{-2}) \right)^{(2t+1/2)/s},
\end{eqnarray*}
for some $C>0$. \\
\\
\textit{Third case}: $a_k \sim k^s$ and $b_k\sim e^{-\gamma k^r}$. Set
$$ D_0 = \left\lfloor \left( \frac{1}{4\gamma} \log \sigma^{-2}) \right)^{1/r} \right\rfloor.$$
Then
\begin{eqnarray*}
\rho_2^2(\mathcal{E}_{a,2}(R),\alpha,\beta)  
& \leq & \rho_{D_0}^2 + R^2 a_{D_0}^{-2},\\
& \leq & \sqrt{D_0} \sigma^2 b_{D_0}^{-2} + R^2 a_{D_0}^{-2}, \\
& \leq & \sigma + C \left(\log(\sigma^{-2}) \right)^{-2s/r} \leq (C+1) \left(\log(\sigma^{-2}) \right)^{-2s/r},
\end{eqnarray*}
for some $C>0$. Concerning the lower bound, we set
$$ D_1 = \left\lceil \left( \frac{1}{2\gamma} \log (\sigma^{-2}) \right)^{1/r} \right\rceil.$$
Then
\begin{eqnarray*}
\rho_2^2(\mathcal{E}_{a,2}(R),\alpha,\beta)
& \geq & \rho_{D_1}^2 \wedge R^2 a_{D_1}^{-2},\\
& \geq & \sigma^2 b_{D_1}^{-2} \wedge R^2 a_{D_1}^{-2}, \\
& \geq & 1 \wedge C \left(\log(\sigma^{-2}) \right)^{-2s/r} = C \left(\log(\sigma^{-2}) \right)^{-2s/r},
\end{eqnarray*}
for some $C>0$. \\
\\
\textit{Fourth case}: $a_k \sim e^{\nu k s}$ and $b_k\sim e^{-\gamma k^r}$. Denote by $\tilde D$ the solution of the equation
$$ \rho_D^2 = R^2 a_D^{-2}.$$
Remark that
$$ \rho_{D_0}^2 \leq R^2 a_{D_0}^{-2} \ \mathrm{where} \ D_0 = \lfloor \tilde D \rfloor,$$
since $(\rho_D^2)_{D\in \mathbb{N}^{\star}}$ and $(a_D^2)_{D\in \mathbb{N}^{\star}}$ are monotone increasing. Hence
$$\rho_2^2(\mathcal{E}_{a,2}(R),\alpha,\beta)  \leq  C \sigma^2 e^{2\gamma D_0^r} + e^{-2\nu D_0^s} \leq (C+1) e^{-2\nu D_0^s}.$$
Then
$$ \rho_{D_1}^2 \geq R^2 a_{D_1}^{-2} \ \mathrm{where} \ D_1 = \lceil \tilde D \rceil.$$
We get:
$$ \rho_2^2(\mathcal{E}_{a,2}(R),\alpha,\beta) \geq \rho_{D_1}^2 \wedge R^2 a_{D_1}^{-2} \geq R^2a_{D_1}^{-2} \geq R^2 e^{-2\nu D_1^s}.$$
In order to conclude the proof, we have to prove that the lower and upper bounds coincide. To this end, remark that $D_1=D_0 +1$. Thus
$$ e^{-2\nu D_1^s} = e^{-2\nu (D_0+1)^s} = e^{-2\nu D_0^s} \times e^{2\nu \lbrace D_0^s - (D_0+1)^s \rbrace} \leq C e^{-2\nu D_0^s},$$
for some constant $C$ as soon as $s\leq 1$.

\subsubsection{Proof of Theorem \ref{thm:lpbod_low}}
The proof will use the one of Theorem \ref{minodure}. We assume that $(\sigma_j)_{j\in J}$ is non-decreasing. Let us first establish a relation between  the $l_p$ ball $\mathcal{E}_{a,p}(R)$ and the sets $\mathcal{S}_{k,n}$. For all $D\in J$, for all $\theta  \in  \mathcal{S}_{ \sqrt{D} , D}$ such that $\|\theta\|_2^2 \leq  \sqrt{D}^{1-2/p} R^2 a_D^{-2}$, we have $\theta \in \mathcal{E}_{a,p}(R)$. Indeed, using H\"{o}lder's inequality
$$ \sum_{j=1}^{+\infty} a_j^p \theta_j^p = \sum_{j:\theta_j \not = 0} a_j^p \theta_j^p \leq (\sqrt{D})^{1-p/2} (\|\theta \|^2)^{p/2} a_D^p \leq R^p.$$
We set $k= \lceil \sqrt{D} \rceil$, $n=D$ and for all $l \in \ac{0,1,\ldots, n-k}$, we define $\theta =(\theta_j,j\in J)$ by (\ref{minoalter}). As pointed out in the proof of Theorem \ref{minodure}, $\theta\in  \mathcal{S}_{k,n}$ and 
$\|\theta\|_2^2 \geq \rho^2$. We also have $\|\theta\|_2^2 \leq \rho^2 \Sigma_{n-k,k}^2/\Sigma_{l,k}^2$. 
This implies that if 
$$ \rho^2 \frac{\Sigma_{n-k,k}^2}{\Sigma_{l,k}^2} \leq  (\sqrt{D})^{1-2/p} R^2 a_D^{-2},$$
 then $\theta\in \mathcal{E}_{a,p}(R)$. \\
Moreover, in the proof of Theorem \ref{minodure}, we proved that if 
$$ \rho^2 \leq  \Sigma_{l,k}^2 \ln\pa{1+\frac{n-l}{k^2}\vee \sqrt{\frac{n-l}{k^2}}}, $$
then
$$\EE_0(L^2_{\mu_{\rho}}(Y)) \leq 1+4(1-\alpha-\beta)^2.$$
This implies by Lemma \ref{lemgene} that $\rho_2^2(\mathcal{E}_{a,p}(R)) \geq \rho^2. $
We finally get
$$\rho_2^2(\mathcal{E}_{a,p}(R)) \geq \Sigma_{l,k}^2 \ln\pa{1+\frac{n-l}{k^2}\vee \sqrt{\frac{n-l}{k^2}}} \wedge
\sqrt{D}^{1-2/p} R^2 a_D^{-2}\frac{\Sigma_{l,k}^2}{\Sigma_{n-k,k}^2}.$$
Since the result holds for all $l \in \ac{0,1,\ldots, n-k}$, we obtain that $\rho_2^2(\mathcal{E}_{a,p}(R)) \geq \rho_1(D)$. 
To obtain that $\rho_2^2(\mathcal{E}_{a,p}(R)) \geq \rho_2(D)$, we consider, as in the proof of Theorem \ref{minodure}, for $k= \lceil \sqrt{D} \rceil$ and  $n=D$
\begin{center}
\begin{tabular}{cll} 
 $\theta_j$&$=\omega_j \sigma_j^2 \rho\pa{\sum_{j=n-k+1}^n \sigma_j^4}^{-1/2} \quad$ &$\forall j\in\{n-k+1,\ldots,n\}$, \\
  & $=0 \quad $&$\forall j\notin\{n-k+1,\ldots,n\}.$
 \end{tabular}
 \end{center}
Since $\rho_2^2(\mathcal{E}_{a,p}(R)) \geq \rho_1(D) \vee\rho_2(D) $ for all $D \in J$, the result follows. 

\subsubsection{Proof of Corollary \ref{corrr}}
For the sake of simplicity, we assume that $\sqrt{D}$ is an integer. 
We derive from Comment 2. of Theorem \ref{minodure} that when $\sigma_j=j^{\gamma}$, then
$$\rho_1(D) \geq C(\gamma) \sqrt{D}^{1-2/p}R^2 a_D^{-2}\wedge \Sigma^2_{D/2,\sqrt{D}}. $$
Moreover,
\beqe
\rho^2_{\sqrt{D},D} &\leq & C \max_{0\leq l\leq D-\sqrt{D}} \Sigma^2_{l,\sqrt{D}}\leq C \sum_{j=D-\sqrt{D}+1}^D j^{2\gamma}\\
&\leq & C D^{2\gamma+1/2} \leq  C \Sigma^2_{D/2,\sqrt{D}}.
\eeqe
When $\sigma_j=\exp^{\gamma j}$, we have from Theorem \ref{thm:lpbod_low} that
$$\rho_2^2(\mathcal{E}_{a,p}(R)) \geq \sup_{D \in J} \rho_2(D) . $$
Moreover, $ \rho_2(D) \geq  \sqrt{D}^{1-2/p}R^2 a_D^{-2}\wedge \sigma_D^2$. We conclude by noticing that, in this case, we also have that $ \rho^2_{\sqrt{D},D} \leq C(\gamma) \sigma_D^2$.

\subsubsection{Proof of Proposition \ref{thm:lpbod_upp}}
It follows from Bonferonis's inequality that $\Phi_{\alpha}^{\dagger}$ is a level-$\alpha$ test. 
Introduce 
$$ A = \left\lbrace  D \in J, R^2 a_D^{-2} \sqrt{D}^{1-p/2} \leq \rho_{\lceil \sqrt{D} \rceil, D}^2 \right\rbrace.$$
In a first time, we suppose that $A$ is empty. From the definition of $D^{\dagger}$, we get $D^{\dagger} = N$ and 
$$ P_{\theta}(\Phi_{\alpha}^{\dagger}=0) \leq P_{\theta}(\Phi_{D^{\dagger},\alpha/2}=0) = P_{\theta}(\Phi_{N,\alpha/2}=0) \leq \beta,$$
for all sequence $\theta$ satisfying
$$ \sum_{j\in J} \theta_j^2 = \| \theta \|^2 \geq C \rho_N^2,$$
for some positive constant $C$. Since $A$ is empty and using the comments following Theorem \ref{minodure}
$$\rho_N^2 \leq C \rho_{\lceil \sqrt{N} \rceil, N}^2  \leq C(\rho_{\lceil \sqrt{N} \rceil, N} ^2 \wedge \sqrt{N}^{1-2/p} a_N^{-2} R^2) \leq C \sup_{D \in J}(\rho_{\lceil \sqrt{D} \rceil, D}^2 \wedge\sqrt{D}^{1-2/p} a_D^{-2} R^2).$$
Hence, our test is powerful as soon as $\|\theta \|^2 \geq \rho^2_{a,p,R}$.\vskip .1in From now on, we assume that the set $A$ is not empty: $D^{\dagger} \leq N$. For all $j\in J$, set
$$ \mu_j = 2(\sqrt{5}+4)\ln \left( \frac{\pi^2 (j- D^{\dagger})^2}{3\alpha \beta} \right).$$
Two different situations may occur:
\begin{description}
\item{1/} For all $j>D^{\dagger}$, $b_j^2\theta_j^2 \leq \sigma^2 \mu_j^2$ for some sequence $\mu_j$, i.e. all the coefficients $\theta_k$ have poor importance after the rank $D^{\dagger}$.
\item{2/} There exists at least $j \in \lbrace D^{\dagger},\dots,N \rbrace$ such that $b_j^2\theta_j^2 \geq \sigma^2 \mu_j^2$, i.e. there exist significant coefficients after the rank $D^{\dagger}$. 
\end{description}
First consider the case 2/. Recall that in this case, the set $A$ is not empty and there exists $j'\in \lbrace D^{\dagger},\dots,N \rbrace$ such that $b_{j'}^2 \theta_{j'}^2 > \sigma^2 \mu_{j'}^2$. In this particular setting, we have to use the threshold test in order to detect these coefficients. More precisely,
$$ P_{\theta}(\Phi_{\alpha}^{\dagger} =0) \leq P_{\theta} \left( \sup_{j>D^{\dagger}} \Phi_{\lbrace j \rbrace, 3\alpha/\pi^2 (j- D^{\dagger})^2} = 0 \right) \leq P_{\theta} \left( \Phi_{\lbrace j' \rbrace, 3\alpha/\pi^2 (j'- D^{\dagger})^2} = 0 \right).$$
Thanks to inequality (29) of \cite{yannick}, we know that this probability is smaller than $\beta$  as soon as:
$$ \theta_{j'}^2 > \sigma^2 b_{j'}^{-2} \ln \left( \frac{\pi^2 (j-D^{\dagger})^2}{3\alpha \beta} \right)2(\sqrt{5}+4).$$
This is exactly the assumption made in case 2/.\vskip .1in Now, we consider point 1/. Let $j> D^{\dagger}$,
\begin{eqnarray*}
\theta_j^2 & = & \theta_j^{2-p} b_j^{2-p} \theta_j^p b_j^{-(2-p)}, \\
& \leq & \sigma^{2-p}\mu_j^{2-p} \theta_j^p b_j^{-(2-p)}.
\end{eqnarray*}
Then, we get
\begin{eqnarray*}
\sum_{j>D^{\dagger}} \theta_j^2 
& \leq & \sigma^{2-p} \sum_{j>D^{\dagger}} \theta_j^p b_j^{-(2-p)} \mu_j^{2-p},\\
& \leq & \sigma^{2-p} \sum_{j>D^{\dagger}} a_j^p \theta_j^p a_j^{-p} b_j^{-(2-p)} \mu_N^{2-p},\\
& \leq & \sigma^{2-p} R^p \max_{j>D^{\dagger}} a_j^{-p} b_j^{-(2-p)} \mu_N^{2-p}.
\end{eqnarray*}
Since the sequence $(a_j^{-p}b_j^{-(2-p)})_{j\in\N}$ is assumed to be monotone non increasing, we can control the bias as follows
$$ \sum_{j>D^{\dagger}} \theta_j^2 \leq \sigma^{2-p} R^p a_{D^{\dagger}}^{-p} b_{D^{\dagger}}^{-(2-p)} \mu_N^{2-p}.$$
In order to conclude the proof, we have to bound the right hand side of the above inequality. First assume that the problem is mildly ill-posed, i.e. $(b_k)_{k\in N^{\star}} \sim (k^{-t})_{k\in\N}$ for some $t>0$. Then
\begin{eqnarray*}
D^{\dagger} & = & \inf \left\lbrace D\in J, R^2 a_D^{-2} (\sqrt{D})^{1-2/p} \leq \sigma^2 D^{2t +1/2} \right\rbrace,\\
& = & \inf \left\lbrace D\in J, R^2 a_D^{-2} \leq \sigma^2 D^{2t +1/p} \right\rbrace.
\end{eqnarray*} 
Thus
\begin{equation}
\sum_{j>D^{\dagger}} \theta_j^2 \leq \sigma^{2-p} \sigma^p (D^{\dagger})^{tp+1/2} (D^{\dagger})^{2t-tp} \mu_N^{2-p} \leq \sigma^2 (D^{\dagger})^{2t+1/2} \mu_N^{2-p}.
\label{eq:prop4_1}
\end{equation}
Hence, we deduce from Proposition \ref{bornesup1} that 
$$ P_f(\Phi_{\alpha}^{\dagger} = 0) \leq P_{\theta}( \Phi_{D^{\dagger},\alpha/2} = 0) \leq \beta,$$
for all sequence $\theta$ satisfying $ \sum_{j=1}^{D^{\dagger}} \theta_j^2 \geq C_{\alpha,\beta} \sigma^2 (D^{\dagger})^{2t+1/2}$, which is equivalent to 
\begin{equation}
\|\theta \|^2 \geq C_{\alpha,\beta} \sigma^2 (D^{\dagger})^{2t+1/2} +  \sum_{j>D^{\dagger}} \theta_j^2.
\label{eq:prop4_2}
\end{equation}
The first point of Proposition \ref{thm:lpbod_upp} follows from (\ref{eq:prop4_1}) and (\ref{eq:prop4_2}). Now assume that the problem is severely ill-posed, i.e. $(b_k)_{k\in\N^{\star}} \sim (e^{-\gamma k})_{k\in\N^{\star}}$ for some positive constant $\gamma$. In this setting,
$$ D^{\dagger} = \inf \left\lbrace D\in J, R^2 a_D^{-2} \sqrt{D}^{1-p/2} \leq \sigma^2 e^{2 \gamma D} \right\rbrace.$$
Hence,
$$ \sum_{j>D^{\dagger}} \theta_j^2 \leq \sigma^2 e^{2 \gamma D} \sqrt{D^{\dagger}}^{1-p/2} \mu_N^{2-p}.$$
An inequality similar to (\ref{eq:prop4_2}) holds, which concludes the second point.

\section{Appendix}

{\bf Proof of Lemma \ref{khideux} : }\\
 We first compute the Laplace transform of $\hT$. 
Easy computations show that for $t<1/(2\sigma_j^2)$,
$$ \EE\cro{\exp(t(\theta_j+\sigma_j \epsilon_j)^2)}= \frac{1}{\sqrt{1-2t\sigma_j^2}}
\exp\pa{\frac{t\theta_j^2}{1-2t\sigma_j^2}}.$$
This implies that for $t<\min_{1\leq j\leq D} 1/(2\sigma_j^2)$,
$$ \EE\cro{\exp(t\hT)}=\exp\pa{\sum_{j=1}^D\frac{t\theta_j^2}{1-2t\sigma_j^2}}\prod_{j=1}^D \frac{1}{\sqrt{1-2t\sigma_j^2}}.$$
Moreover,
$$ \EE(\hT)=\sum_{j=1}^D \theta_j^2+ \sum_{j=1}^D \sigma_j^2.$$
This leads to
\beqe
\EE\cro{\exp(t(\hT- \EE(\hT))}&= &\exp\pa{\sum_{j=1}^D\frac{2t^2\theta_j^2 \sigma_j^2}{1-2t\sigma_j^2} -t\sigma_j^2}\prod_{j=1}^D \frac{1}{\sqrt{1-2t\sigma_j^2}}\\
&=&  \exp\pa{\sum_{j=1}^D\frac{2t^2\theta_j^2 \sigma_j^2}{1-2t\sigma_j^2} -t\sigma_j^2}\exp\pa{-\frac{1}{2}\sum_{j=1}^D \log(1-2t\sigma_j^2)}.
\eeqe
We use the following inequality which holds for $x <1/2$ :
\beq \label{Ineglog}
x \cro{\frac{1}{2}\log \pa{1-2x} +x +\frac{x^2}{1-2x} }\geq 0.
\eeq
This inequality implies that for all $t < \min_{1\leq j\leq D} 1/2\sigma_j^2$, 
$$ \log \EE\cro{\exp(t(\hT- \EE(\hT))}\leq \sum_{j=1}^D \frac{t^2\sigma_j^4}{1-2t\sigma_j^2}+2t^2 \sum_{j=1}^D \frac{\theta_j^2 \sigma_j^2}{1-2t\sigma_j^2}.$$
This leads to 
$$\log \EE\cro{\exp(t(\hT- \EE(\hT))}\leq  \frac{ t^2 \Sigma}{1-2t \sup_{1\leq j \leq D} ( \sigma_j^2)}.$$
We now use the following lemma which is proved in Birg\'e \cite{Lucien} (see Lemma 8.2) : 
\begin{lemma}
Let $X$ be a random variable such that 
$$
\log\left(\EE\cro{\exp(tX)}\right)\leq\frac{(at)^2}{1-bt}\quad\mbox{for }
0<t<1/b$$
where $a$ and $b$ are positive constants. Then
$$\PP\pa{X\ge 2a\sqrt{x}+bx}\leq\exp(-x)\quad\mbox{for all }x>0.$$
\end{lemma}
Hence, inequality (\ref{majo}) is proved. Let us now prove inequality (\ref{mino}). \\
For all $z \in \mathbb{R}$, 
\beqe
\PP\pa{\hT -\EE(\hT)\leq - z}&= & \PP\pa{ -\hT+\EE(\hT) -z \geq 0}\\
&\leq& \inf_{t>0} \EE\pa{ e^{t(-\hT+\EE(\hT) -z )}}\\
&\leq & \inf_{t<0} \EE\pa{ e^{t(\hT-\EE(\hT) +z )}}.
\eeqe
We have, from the above computations 
$$ \ln\pa{\EE\pa{ e^{t(\hT-\EE(\hT) +z )}}}  = \sum_{j=1}^D \cro{\frac{2t^2 \theta_j^2 \sigma_j^2}{1-2t\sigma_j^2} -t \sigma_j^2 -\frac{1}{2} \ln(1-2t \sigma_j^2)}+t z.$$
We now use (\ref{Ineglog}) for $x= t \sigma_j^2$ with $t<0$. We obtain 
$$ \frac{1}{2}\ln (1-2t \sigma_j^2)+ t \sigma_j^2 + \frac{t^2 \sigma_j^4}{1-2t \sigma_j^2} \leq 0.$$
This implies that
$$ \frac{2t^2 \theta_j^2 \sigma_j^2}{1-2t\sigma_j^2} \leq -2t \theta_j^2 -\frac{\theta_j^2}{\sigma_j^2}\ln(1-2t \sigma_j^2).$$
Hence, for all $t<0$, $z \in \mathbb{R}$, 
$$\EE\pa{ e^{t(\hT-\EE(\hT) +z )}} \leq  \exp\cro{-\sum_{j=1}^D \pa{\frac{1}{2}\log\pa{1-2t\sigma_j^2}+t\sigma_j^2}\pa{1+2\frac{\theta_j^2}{\sigma_j^2}}+tz}.$$
We use this inequality with  $z= 2 \sqrt {\Sigma x}$, and $t_x=-\sqrt{x}/\sqrt{\Sigma}$. 
$$ \PP\pa{\hT -\EE(\hT)\leq -2 \sqrt {\Sigma x}} \leq \EE\pa{ e^{t_x(\hT-\EE(\hT) +  2 \sqrt {\Sigma x})}}.$$
Moreover, 
$$\EE\pa{ e^{t_x(\hT-\EE(\hT) +  2 \sqrt {\Sigma x})}} = \exp\cro{-\sum_{j=1}^D \pa{\frac{1}{2}\log\pa{1-2 \frac{\sqrt{x}}{\sqrt{\Sigma}}\sigma_j^2} -\frac{\sqrt{x}}{\sqrt{\Sigma}} \sigma_j^2}\pa{1+2\frac{\theta_j^2}{\sigma_j^2}}-2x }.$$
We use the following inequality which holds for all $u \geq 0$ : 
$$ \frac{1}{2} \log(1+ 2u) -u \geq -u^2,$$
and we apply this  inequality to $ u= - t_x \sigma_j^2$. We obtain that for all $x\geq 0$, 
$$ \PP\pa{\hT -\EE(\hT)\leq - 2 \sqrt {\Sigma x}} \leq \exp(-x).$$
This concludes the proof of Lemma \ref{khideux}.


\begin{thebibliography}{17}
\bibitem{Aldous} Aldous, D. J. (1985) {\it Exchangeability and related topics.} Ecole d'\'et\'e de probabilit\'es de Saint-Flour XIII, Lect. Notes Math. 11117, 1-198. 
\bibitem{yannick} Baraud, Y. (2002) {\it Non asymptotic minimax rates of testing in signal detection}, Bernoulli, {\bf 8}, 577-606.
\bibitem{BHL} Baraud, Y., Huet, S., and Laurent, B. (2003) {\it Adaptive tests of linear hypotheses by model selection}, Ann. Statist., {\bf 31}, no. 1, 225-251.
\bibitem{Lucien} Birg\'e, L.  (2001) {\it An alternative point of view on Lepski's method}, State of the art in probability and statistics (Leiden, 1999) (ed. Monogr., IMS Lecture Notes, {\bf 36}, 113-133.
\bibitem{munkpreprint}Bissantz, N. and Claeskens, N and Holzmann, H. and Munk, A. (2008) {\it Testing for lack of fit in inverse regression, with applications to
  biophotonic imaging}, { preprint}.
\bibitem{Butucea} Butucea, C. (2007) {\it Goodness-of-fit testing and quadratic functional estimation from indirect observations}, Ann. Statist., {\bf 35}, no. 5, 1907-1930.
\bibitem{MR2421941} Cavalier, L. (2008) {\it Nonparametric statistical inverse problems}, Inverse Problems, 24(3).  
\bibitem{Donoho} Donoho, D. (1995), {\it Nonlinear solution of linear inverse problems by Wavelet-Vaguelette decomposition}, Applied and computational harmonic analysis, {\bf 2}, 101-126.
\bibitem{FLR} Fromont, M., Laurent, B. and Reynaud-Bouret, P.  (2009) {\it Adaptive test of homogeneity for a Poisson process }, ArXiv:0905.0989v1, Submitted. 
\bibitem{Erm} Ermakov, M. S. (2006), {\it Minimax detection of a signal in the heteroscedastic Gaussian white noise}, J. Math. Sci., {\bf 137}, No. 1, 4516-4524. 
\bibitem{Ing} Ingster, Yu.I. (1993)  {\it Asymptotically minimax testing for nonparametric alternatives I-II-III}, Math. Methods Statist., {\bf 2}, 85-114, 171-189, 249-268.
\bibitem{IngSap} Ingster, Yu.I., Sapatinas, T. and Suslina, I.A. (2010) {\it Minimax signal detection in ill-posed inverse problems}. Working paper.
\bibitem{IngSuslina} Ingster, Yu.I. and Suslina, I.A. (1998) {\it Minimax detection of a signal for Besov bodies and balls}, Problems Inform. Transmission, {\bf 34}, 48-59. 
\bibitem{IngBook} Ingster, Yu.I. and Suslina, I.A. (2002) {\it Nonparametric goodness-of-fit testing under Gaussian models}, Lecture Notes in Statistics, 169. Springer-Verlag, New York. 
\bibitem{BLPM} Laurent, B. and Massart, P. (2000) {\it Adaptive estimation of a quadratic functional by model selection}, Ann. Statist., {\bf 28}, no. 5, 1302-1338.
\bibitem{autre} Laurent, B. and Loubes, J-M. and Marteau, C. (2010) {\it Testing inverse problems: a direct or an indirect problem}, preprint.
\bibitem{Lepspok} Lepski, O. V., and Spokoiny, V. G. (1999) {\it Minimax nonparametric hypothesis testing: the case of inhomogeneous alternative}, Bernoulli, {\bf 5}, 333-358. 
\bibitem{loublud1} J-M. Loubes and C.~Ludena. (2008) {\it Adaptive complexity regularization for inverse problems},
{Electronic Journal Of Statistics}, 2:661--677.
\bibitem{osulliv}O'Sullivan, F. (1986) {\it  A statistical perspective on ill-posed inverse problems}, { Statist. Sci.}, 1({\bf 4}), 502--527.
\end{thebibliography}
\end{document}